\documentclass[a4paper,10pt]{article}
\usepackage{amsmath}
\usepackage{amsfonts}
\usepackage{graphicx}
\usepackage{dsfont}
\usepackage[english]{babel}
\usepackage[ansinew]{inputenc}
\usepackage{verbatim}
\usepackage{afterpage}
\usepackage{longtable}
\usepackage{amssymb}

\newcommand{\be}{\begin{equation}}
\newcommand{\ee}{\end{equation}}
\newcommand{\bea}{\begin{equation*}}
\newcommand{\eea}{\end{equation*}}

\newcommand{\dec}{\mathcal D}

\newtheorem{theorem}{Theorem}
\newtheorem{definition}[theorem]{Definition}
\newtheorem{proposition}[theorem]{Proposition}

\newtheorem{remark}{Remark}
\newtheorem{ass}{Assumption}
\def\proof{\noindent{\bf Proof}\ \ }
\def\qed{\hfill \vrule height 7pt width 7pt depth 0pt
		\medskip}

\def\N{\mathds{N}}
\def\R{\mathds{R}}

\def\MM{\mathrm M} 
\def\MA{\mathrm A} 
\def\MB{\mathrm B} 
\def\MC {\mathrm C} 
\def\zz{\zeta_0}  
\def\zu{\zeta_1} 
\def\dd{\mathrm{d}}

\def\erfc{\text{erfc}}
\def\SNR{\text{SNR}}
\def\DEP{\text{DEP}}
\def\EDP{\text{EDP}}
\title{A Decoding Approach to Fault Tolerant Control of Linear Systems with Quantized Disturbance Input}
\author{Sophie M. Fosson}

\begin{document}
\maketitle 
\abstractname{.}
The aim of this paper is to propose an alternative method to solve a Fault Tolerant Control problem. The model is a linear
system affected by a disturbance term: this represents a large class of technological faulty processes. The goal is to make the system
able
to tolerate the undesired perturbation, i.e., to remove or at least reduce its negative effects; such a task is performed in three
steps: the
detection of the fault, its identification and the consequent process recovery. When the disturbance function is known to be
\emph{quantized} over a
finite number of levels, the detection can be successfully executed by a recursive \emph{decoding} algorithm, arising from Information
and Coding Theory and suitably adapted to the control framework. This technique is analyzed and tested in a flight control
issue; both theoretical considerations and simulations are reported.

\section{Introduction}
Fault Tolerant Control (FTC for short, \cite{bl:06},\cite{is:06},\cite{du:09}) aims to cancel or contain the consequences of faults in
an automation
system. Such an operation is fundamental in modern technological processes, which are required to assure
robust performance, stability and safety even in case of partial malfunctions or degradations. Often, robustness is
achieved by redundancy, say by the introduction of many control components like sensors; nevertheless, this sophistication 
naturally increases the probability of breakdown and then continues to motivate
the research on reliable control systems.\\

The problem of upholding the functionality of an apparatus affected a disturbance is ubiquitous in the industrial and transport
fields. In particular, FTC systems are widely applied in those contexts where human health and environment are concerned, for example, 
in the design of mechanical and chemical  plants; nuclear power reactors; medical systems; aircrafts, helicopters and spacecrafts; 
automotive engines, railway and marine vehicles. Another interesting application is
in the communication networks (for instance, wireless sensor networks),  where the aim of FTC is to avoid unexpected interruptions of
data flow in case of troubled connectivity or impaired nodes. In all these contexts, a satisfying FTC design can prevent
non-reversible failures and stops, with the ultimate objective of reducing health, environmental and economic damages.

The literature about FTC is definetely widespread and contributions arise from diverse applied mathematical domains. 
In order to get into the argument, there are many survey works that
introduce the main theoretical concepts and
provide classifications of the outstanding FTC approaches, with detailed references. For example, we refer the reader to the
recent
review \cite{bib_rev}, which supplies a comprehensive bibliography, and to \cite{ji:05}, \cite{st:91}, \cite{pat:97}, \cite{wil:76}.

As far as the applications are concerned, aircraft flight control has been motivating FTC research since 1970s, given the evident
danger that aircraft faults may cause to human safety.  Therefore, a significant amount of papers has been produced on the argument,
 taking account of the wide variety of issues and  models introduced in the study of flight dynamics. For a general  overview see
\cite{st:05}, \cite{et:85} and the up-to-date book \cite{du:09} that in Chapter II provides the list
of the most
common flight control systems, with the relative references. \\

In this work, a linear model with a multiplicative disturbance factor is considered, which is very
common in flight framework (\cite{yy:06}); in particular, we will adopt a system presented in  \cite{ac:85},\cite{ac:84}
and studied also in \cite{yee:02}, \cite{fag:04} as an application test.\\

Even if FTC systems can be designed in many different ways according to the specific aim they are conceived for,  in general
they all have
to perform the following main tasks:
\begin{enumerate}
 \item the Fault Detection, i.e., the controller makes a binary decision on the presence of a malfunction;
\item the Fault Identification, i.e., the controller determines or  estimates the size of the disturbance; if necessary, Identification
is preceded by Fault Isolation, that is, the location of the impaired component;
\item the eventual active compensation to the fault, i.e., the reconfiguration of the system inputs and/or parameters in order to
maintain, as much as possible, the integrity of the process.
\end{enumerate}
Fault Detection and Identification (FDI) can be undertaken in diverse ways. In the cited works, in particular \cite{du:09} a
comprehensive discussion about the most popular FDI schemes is presented: among them, we remind the unknown input observers (UIO,
\cite{PatChen}, \cite{Vis87}) and
residual generation, Kalman filtering, the statistical methods and the more recent techniques based on neural networks (\cite{Nap95}).

This paper is devoted to the case when a quantized disturbance input is introduced in a continuous linear system. Such an \emph{hybrid}
model, which combines discrete and continuous dynamics, is motivated by the
upcoming digitalization of modern devices: a quantized disturbance may represent the
switches of  actuators or sensors and the malfunctions in digital components; moreover, it may describe the behavior of any
mechanical device that is known to occupy only certain positions and also the approximation of a continuous disturbance. 

Results about FTC for  hybrid systems are not very common.
In part, they
can be retrieved in the extensive discussion about the detection of \emph{abrupt} changes in dynamical
systems, whose leading work is \cite{bas} (while some further contributions are given by \cite{abrupt:1999} and \cite{abrupt:2000}). The
problem of estimating brusque alterations is always actual (as an example, see
\cite{abr_radar} and \cite{abr_medicalimaging}, which respectively concern medical imaging and ground-penetrating radar issues) and
in general is approched by classical estimation techniques, such as Kalman Filtering.\\
Recently, input quantization in linear systems has been studied  in particular with the aim of reducing the effects of a coarse
quantization (\cite{park}, \cite{eliamitter}). In this work, instead, our purpose  is exploiting the information that the disturbance
input is quantized to detect the fault occurrene: it follows that quantization is supposed to be already performed in a satisfactory
way.\\

In order to evaluate the quantized input disturbance, an original Information theoretic approach is proposed in this paper: given the
discrete nature of the disturbance, FDI is
performed by a \emph{decoding technique} derived from the framework of digital transmissions and Coding Theory (\cite{RU}).
The algorithm we will introduce has already been tested in Deconvolution issues (\cite{09}). The problem we address here still is a
Deconvolution problem,
given that we assume a linear system as model, but in addition a compensation task is introduced to minimize the consequence of faults:
our FTC is conceived with a feedback loop that supplies a compensation input in real-time and then continuously reconfigures the system
(which naturally does not happen in classical Deconvolution
issues).

The structure of
the paper is the following: in Section II, we describe the problem we aim to study; in Section III, we introduce the decoding algorithm
furtherly used for the Fault Detection; in Section IV, we provide a theoretical analysis of the algorithm in terms of minimization of a
suitably defined \emph{Error Function} that represents the distance between the optimal
behavior (i.e., without disturbance) and the output of the FTC itself;  sensitivity to the false alarm (\emph{false positive}) and to
miss fault detection (\emph{false negative});  promptness of detection and reconfiguration. In Section V, wi give the design criteria
to  obtain the best performance from our algorithm, while in Section VI we show a few significant
simulations about a specific numerical example, arisen from Flight Control literature; finally,  Section VII is devoted to some
conclusive observations.
\subsection{Notation}
In this paper, the following notation will be used:
\begin{itemize}
\item given a subset $A$ of a set $X$, $\mathds{1}_A:X\to\{0,1\}$ will denote the indicator function, defined by $\mathds{1}_A(x)=1$ if
$x$ belongs to $A$ and $\mathds{1}_A(x)=0$ otherwise;
\item the function $\erfc$ is defined by  $\erfc(x)=\int_x^{+\infty}e^{-s}ds$ for any $x \in \R$;
\item random variables will be indicated by capital letters;
\item given any variable $X$, $\hat{x}$ will denote its estimation.
\end{itemize}

  \section{Problem Statement}
  In this paper, we consider processes that can be modeled by the following linear, finite-dimensional system:
    \be\label{fault}\left\{\begin{array}{l} \dot{x}(t)=\MA x(t)+\MB z(t)f(t)~~~t\in[0,T]\\
    x(0)=0\\
    y(t)=\MC x(t)\end{array}~~~~~\right.\ee
    where $x(t) \in \R^n$, $y(t) \in \R^m$, $f(t)$ and $z(t)$ are scalar functions and $\MA$, $\MB $ and $\MC$ are constant matrices
with
  consistent dimensions. $f(t)$ is a known input signal, while $z(t)$ is a disturbance  modelling some fault in the system. Typically,
$z(t)\in(0,1]$; if
  $z(t)=1$,  the system operates in its nominal regime and is totally driven by $f(t)$: this is the condition that one aims to
  reproduce even when $z(t) \in (0, 1)$, i.e., when some unexpected breakdown, interruption or loss of effectiveness affects the
dynamics.

  In order to achieve that, a control input $u$ is introduced, which adjusts the dynamics as follows:
    \be\label{fault2} \dot{x}(t)=\MA x(t)+\MB z(t)\left(f(t)+u(t)\right)\ee

  Notice that to maintain the error-free behavior, say $\MB z(t)\left(f(t)+u(t)\right)=\MB f(t)$, in principle it is
  sufficient to fix
  $u(t)=f(t)\left(\frac{1}{z(t)}-1\right)$, but, in the real applications, this is often impossible for the following motivations.
 Generally,  the disturbance $z$ is not known and the the controller can access it only through the observation of the
  output $y$. In order to determine $z$ one has to perform a \emph{deconvolution}, that is, to invert the solution of equation 
  (\ref{fault2}) with initial condition $x(0)=0$:
  \be\label{fault2_sol} y(t)=\MC x(t)=\MC \int_0^t e^{(t-s)\MA }\MB z(s)(f(s)+u(s)) ds\ee
Furthermore, the acquisition of the data usually is not exact. This inaccuracy can be modeled by an additive noise $n(t)$ in
the output (in this work, $n(t)$ will be defined as a white gaussian
noise): the available function now is $r(t)=y(t)+n(t)$. 

Under this condition, the inversion of expression (\ref{fault2_sol}) becomes tricky: deconvolution is in fact known to be an ill-posed
and ill-conditioned problem, that is, the uniqueness of solution is not guaranteed and also small errors in the data may raise large
errors in the solution. In conclusion, the reconstruction of $z(t)$ by inversion may produce outcomes very far from the correct ones;
for this reason, an estimation approach to the problem is the most suitable one.
 
In addition to that, in this work we make the following
  The controller can access $y$ only at certain time instants, say each $\tau$ time instants. Hence, the available data are the
samples $r_k=r(k\tau)$ where $K\in\N$, $k\in\{0,\dots,K-1\}$ (for simplicity, let us suppose that $K\tau =T$).
 
Moreover, in this work, two further main assumptions are made.
\begin{ass}
The controller can access $r(t)$ only at  each $\tau$ time instants. The available data are 
the samples $r_k=r(k\tau)$ where $k\in\{1,\dots,K\}$ and $K\in\N$ is supposed to be such that $K\tau =T$.
\end{ass}
\begin{ass}
The disturbance function $z(t)$ is known to be quantized over two levels, say $z(t)$ can assume only two values $\zz$ and $\zu$.
\end{ass}

$\zz$ and $\zu$ may respectively represent the nominal and the faulty conditions ($\zz=1$, $\zu\in(0,1)$). Such a binary situation
naturally occurs in many engineering applications: it can model, for instance, the abrupt
blocking of an actuator, the sharp loss of efficiency of a device, the sudden disconnection of some component, the functioning of alarm
sensors. In the next, we will generally refer to the jumps from $\zz$ and $\zu$ and vice-versa as \emph{switch points}.\\
Notice that Fault Detection and Identification are coincident under this assumption: the decision on the fault presence automatically
determines also its size. \\

In this work, being aware of all these conditions, we aim to estimate $z(t)$  as well as possible in order to provide the best control
input to the system. Clearly, the estimation has to be performed on-line, that is, each time a sample is acquired (notice that the
sampling inevitably undertakes some delay): each $\tau$ instants the controller tries to detect eventual faults and
consequently updates the system design.

For mathematical simplicity, the eventual switch points of $z(t)$ are suppoed to occur at the time instants
$k\tau$, in order to have synchronization with the output sampling. Hence, we can write:
  \be\label{quantiz} z(t)=\sum_{k=0}^{K-1} z_k \mathds{1}_{[k\tau,(k+1)\tau[}(t)~~~z_k \in \left\{\zz,\zu\right\}\ee
Now, $z(t)$ is equivalent to the binary sequence $(z_0,\dots,z_{K-1})\in\{\zz,\zu\}^K$: the estimation problem is actually discrete.
Let $\hat{z}_k$ be an estimate
 of $z_k$: since the operation must be performed on-line, we expect $\hat{z}_{k-1}=\mathcal{D}(r_1,\dots,r_{k})$,
  where $\mathcal{D}$ indicates a detection/estimation function.

  Taking account of the conditions mentioned before, the natural definition of the control input is:
 \be
  u(t)=f(t)\left(\frac{1}{\hat{z}_{k-1}}-1\right)\mathds{1}_{[k\tau,(k+1)\tau)}(t)~~~~~k=0,\dots,K-1\ee 
  $u(t)$ is computed and introduced in the system each $\tau$ time instants. Consider now a generic interval $[k\tau,(k+1)\tau)$. Being
based on the estimate $\hat{z}_{k-1}$ relative to the previous interval, $u(t)$ is deceptive when a switch occurs at $k\tau$: the delay
$\tau$ underlies a temporary, unavoidable deviation (even in case of correct detection) from the right trajectory. This issue will be
widely discussed in
  the next; for the moment, let us just observe that switch points cause the most of the problems in our FTC model. For this reason,
  \emph{permanent} interruptions, i.e., \emph{failures} (which involve just one switch point) are definitely preferable than
  \emph{transient} faults for our purpose, though this should appear as a paradox in the practice.

\subsection{Illustrative Example: a Flight Control Problem}\label{example}
A typical example of FTC problem arises from the literature of Flight Control. Systems of kind (\ref{fault}) are often used to model
different aspects of the aerospace dynamics. For instance, if we consider the matrices
  \be \MA =\left[\begin{array}{ccc}-0.5162& 26.96&178.9\\-0.6896&-1.225&-30.38\\0&0&-14 \end{array}\right] \ee
    \be \MB =\left[\begin{array}{c}-175.6\\0\\14\end{array}\right],~~~\MC =\left[1~~ 12.43~~ 0\right]\ee
  the system (\ref{fault}) represents the longitudinal short-period mode of an F4-E jet with additional horizontal canards, in
supersonic conditions. The vector $x$ determines the longitudinal trajectory: its three entries respectively
represent the normal
acceleration, the
pitch rate and the deviation of elevator deflection from the trim position. The output
  $y(t)$ is the $C^*$ response, a usual parameter in flight mechanics that synthesizes the aircraft response to the pilot inputs;
typically, the $C^*$ response must lie in a given admissible envelope.\

This application example is illustrated in the Appendix D.1 of \cite{ac:85} and studied also in
\cite{ac:84},\cite{yee:02},\cite{fag:04}.

In this context, $f(t)$ can be interpreted as the elevator deflection command and $z(t)$  as the indicator of
the status of the elevators: $z=\zz$ may attest a good status, while the switch to $z=\zu$ may denote an abrupt loss of
effectiveness. In such a case, the controller is required to detect the accident and add the suitable control input $u(t)$ in order to
recover the optimal trajectory, say the one imposed by the flight plan. In terms of the output $y(t)$, one aims to maintain or to
bringit it back into the prescribed envelope.

Notice that in this case, it makes sense to suppose the fault to be definitive, that is, the elevator cannot recover its
efficiency during the flight. We then talk about a failure. This situation often occurs in the applications, which motivates us
to
focus on it in our following analysis. This Flight Problem will be retrieved later and used as test application for the implementation
of our detection algorithm, which is
introduced in the next section.

  \section{Fault Detection: The One State Algorithm}
  Given the quantization of  $z_k \in \{\zz,\zu\}$, it makes sense to settle the same set for the estimation: $\hat{z}_k \in
  \{\zz,\zu\}$. This consideration arises from coding/decoding techniques in digital transmissions, where unknown input messages, that
are
  combinations of symbols from a known finite alphabet, must be recovered  within the same alphabet. In other terms, the \emph{decoder}
  is an estimator that exploits the prior information about the input source.

 The detection method that we introduce in this section is derived from an optimal decoding algorithm named BCJR
 after its authors Bahl, Cocke, Jelinek and Raviv (see
\cite{BCJR}). Given the noisy output of a digital transmission, the BCJR computes the probabilities of all the possible codewords,
implementing a maximum a posteriori (MAP, \cite{RU}) estimation through a recursive procedure. In particular, given codes defined on
trellises, it
evaluates the a posteriori probabilities of each state. 

The classical version of the algorithm is constituted by two recursions (one forward, one backward) and requires the transmission of the
whole message   before decoding. Moreover, it also requires the system to have a finite number of states.
Nevertheless, it is possibile to modify the proceudre to avoid these bonds: in spite of reliability, one can make it causal (hence to
work on line) by removing the backward recursion and also it can be simplified by considering not all the possible states, but just a
fixed number of states. In \cite{09}, these variations are widely discussed. The algorithm we introduce here is exactly a causal BCJR
considering just one state at each step (for this reason we refer to it as the One
State Algorithm). The compuations of the probabilities is in this case straightforward and reduces to the comparison between two
Euclidean distances at each step. This makes the algorithm definetely low-complexity, which encourages its implementation. Its
performance actually depends on the specific application case and will be analysed in the next sections. 

Now, let us describe the operative structur of the One State Agorithm in detail.
\subsection{One State Algorithm's pattern}
  Before showing the algorithm, notice that the solution of the equation (\ref{fault2_sol}) can be
  written recursively as
  \be\begin{split} x_k&=e^{\tau \MA }x_{k-1}+z_{k-1}(1-u_{k-1}) \int_0^{\tau}e^{s\MA }\MB  f(k\tau -s)ds\\
  &=e^{\tau \MA }x_{k-1}+\frac{z_{k-1}}{\hat{z}_{k-2}} \int_0^{\tau}e^{s\MA }\MB  f(k\tau-s)ds\\
  x_0&=0\end{split}\ee
 where $x_k=x(k\tau)$, $k=0,\dots,K$. Now, the key idea of the One State procedure is to provide a recursive estimation of the
state $x_k$ and of $z_{k-1}$ given the current lecture $r_k$ and  the estimate of the previous state
  $x_{k-1}$.

  In the next, let us use the following notation: $n_k=n(k\tau)$, $d_E$ indicates the Euclidean
  distance and finally:
\be \MM_{\tau,k}=\int_0^{\tau}e^{s\MA }\MB f(k\tau -s)ds\ee

The One State Algorithm's pattern is then the following:
    \begin{enumerate}
	\item $k=0$. Initialization:  $\hat{x}_0=0$;
	\item $k=1$.\\
  System evolution (with no compensation): $x_1=z_0 \MM_{\tau,1}$.\\
  Lecture: $r_1=y_1+n_1=\MC x_1+n_1$.\\
  Disturbance Estimation: $\hat{z}_0=\left\{\begin{array}{ll}
    \zz &\text{if  }  d_E(r_1,\zz\MC \MM_{\tau,1})\leq d_E(r_1,\zu\MC  \MM_{\tau,1})  \\
    \zu &\text{otherwise  }\\
    \end{array}\right.$\\
    State Estimation : $\hat{x}_1=\hat{z}_0 \MM_{\tau,1}$.

    \item $k=2,\dots,K$.\\
  System evolution (with compensation): $x_k=e^{\tau \MA }x_{k-1}+\frac{z_{k-1}}{\hat{z}_{k-2}}
  \MM_{\tau,k}$.\\
  Lecture: $r_k=y_k+n_k=\MC x_k+n_k$.\\
  Disturbance Estimation: 
  $\hat{z}_{k-1}=\left\{\begin{array}{ll}
  \zz &\text{if  }  d_E(r_k,\MC e^{\tau \MA }\hat{x}_{k-1}+\frac{\zz}{\hat{z}_{k-2}}\MC \MM_{\tau,k})\\&\leq d_E(r_k,\MC e^{\tau
  \MA }\hat{x}_{k-1}+\frac{\zu}{\hat{z}_{k-2}}\MC \MM_{\tau,k}      ) \\
    z_1&\text{otherwise  }\\
    \end{array}\right.$\\
  State Estimation: $ \hat{x}_k=e^{\tau \MA }\hat{x}_{k-1}+\frac{\hat{z}_{k-1}}{{\hat{z}_{k-2}}}\MM_{\tau,k}$.
  \end{enumerate}
  Notice that the system does not have compensation in the first interval $[0,\tau)$, as the first useful lecture is
  performed at time $t=\tau$. 
  For the binary nature of each $z_k$, the process of estimation/detection reduces here to the comparison of two distances.
Moreover,  the storage required is of two locations (one float for the current state and one boolean for the current disturbance): the
algorithm is  definitely low-complexity.

  \section{Theoretical Analysis of the  One State Algorithm}\label{theoanalysis}

This section is devoted to the theoretical description of the behavior and performance of the One State Algorithm applied to the system (\ref{fault})-(\ref{quantiz}) with a failure, that is, there exists a time instant $T_F=k_F\tau\in [0,T]$, $k_F\in \N$ such that
\be\label{failure} z(t)=\left\{\begin{array}{ll}
                 \zz=1 & t\in[0, T_F)\\
                 \zu\in (0,1) & t\in [T_F,T]\\
                \end{array}\right. \ee
or equivalently, $z_k=\zz$ for $k=0,1,\dots,k_F-1$ and $z_k=\zu$ for $k=k_F,1,\dots,K-1$. Switch points are particulary tricky and the
choice to focus on a system with just one switch point allows to isolate the problem and to understand completely the consequences of a
switch. On the other hand, this case is crucial for the applications, where the problem of failures is dramatically serious.

 Our model can be naturally described in probabilistic terms:  the fact that
  lecture noise is supposed to be white gaussian, (that is, a sequence of independent gaussian random variables $N_k\sim\mathcal
  N(0,\sigma^2)$) introduces some amount of uncertainty in the system. In particular, also $\hat{z}$, $x$, $y$, $r$, $\hat{x}$ are
random variables, as they are directly or indirectly functions of
the
noise. To emphasize that stochastic nature, from now
onwards, we will indicate random variables by capital letters. Let us resume the complete recursive system in probabilistic terms:
\be\begin{split}
&X_0=0,~\hat{X}_0=0,~\hat{Z}_{-1}=\zz=1\\
&X_k=e^{\tau \MA }X_{k-1}+\frac{z_{k-1}}{\hat{Z}_{k-2}}\MM_{\tau,k}\\ 
&Y_k=\MC X_k\\
&R_k=Y_k+N_K\\
&\hat{Z}_{k-1}=\dec_1( R_k,\hat{X}_{k-1},\hat{Z}_{k-2})\\
&\hat{X}_k=e^{\tau \MA }\hat{X}_{k-1}+\frac{\hat{Z}_{k-1}}{\hat{Z}_{k-2}}\MM_{\tau,k}, ~~k=1,\dots,K\\ 
\end{split}\ee
where  $\dec_1$  indicates the One State detection function. Notice that $X_0$, $\hat{Z}_{-1}$, $X_1$,
$Y_1$ are actually deterministic, in particular, fixing $\hat{Z}_{-1}=\zz=1$ is just an other way to state that there is no
compensation for the system in the first interval $[0,\tau)$.\\

Finally, we remark that $z(t)$ is not supposed to be driven by some probabilistic law. Such an information
on the input might be useful to improve the detection and has been studied in other deconvolutio contexts (see, for instance,
\cite{09}). Nevertheless, in this work we rather prefer to focus on a specific disturbance.

\subsection{The Error Function}\label{EF}
The performace of the algorithm must be determined through the evaluation of a suitable \emph{error function}, say a distance between
the desired and the real
trajectories. In this work, we adopt as error function the discrete stochastic process $(E_k)_{k=0,1,\dots}$ that describes
the signed distance between the
trajectory of the system with control and compensation $X_k$ and the nominal trajectory $x^N(t)$, at time instants $k\tau$,
$k=0,1,\dots$:
 \be\label{defErr} \left\{\begin{array}{l}
      E_{k}=X_k-x^N(k\tau)\\
     ~~~~~ =e^{\tau \MA }E_{k-1}+\left(\frac{z_{k-1}}{\hat{Z}_{k-2}}-1\right)\MM_{\tau,k}~~~k=1,\dots,K\\
     E_0=0.\\
      \end{array}\right.\ee

The so-defined error function is characterized by the following fact:.
\begin{proposition}\label{equiv_events}
For any $k_0,n\in \N$,  the event $\{E_{k_0+n}=e^{n\tau A}E_{k_0}\}$ corresponds to the event $\{\hat{Z}_{k-1}=z_{k}$ for all
$k=k_0,k_0+1,\dots k_0+n]\}$.
\end{proposition}
\proof It immediately follows from the definition of $E_k$:  for any $n\in\N$, 
the event  $\{E_{k+1}=e^{\tau \MA
}E_{k}\}$ is equivalent to   $\{\hat{Z}_{k-1}=z_{k}\}$ and then $\{E_{k_0+n}=e^{n\tau \MA
}E_{k_0}\}$ corresponds to the event   $\{\hat{Z}_{k_0-1}=z_{k_0}, \hat{Z}_{k_0}= z_{k_0+1},\dots,\hat{Z}_{k_0+n-1}=z_{k_0+n}\}$. \qed\\
Notice that under the hypothesis of the proposition and if $\MA$ is asymptotically stable, $E_k$ exponentially decades to
zero, regardless of the initial value $E_{k_0}$.\\
Moreover, observe that the condition $\hat{Z}_{k-1}=z_{k}$ is not the event of correct detection $\hat{Z}_{k}=z_{k}$, since the feedback
in the system implies a delay $\tau$; however, if $z_k$ is constant over the considered interval, the two events are the same. In the
next, we will focus on this context of constant disturbance, which models the state of the system before and after an irreversible
failure. In particular, we will study the conditions to obtain correct detection, which leads to the exponential decay of the error; 
we will show that even if we cannot achieve the certainty of decodig exaclty in the presence of noise, however we can approximate this
condition satisfactorily, that is, with a probability close to one, at least in some common situations.

More precisely, our goal is to study the probability of the event $E_{k_0+n}=e^{n\tau \MA }E_{k_0}$ conditioned to the fact that
$z_k \text{ constant for any} k\in[k_0,k_0+n]$ and given some  initial conditions at $k_0$ concerning the state of the algorithm, which
will be defined later. In particular, we will find out the conditions that make this probability sufficiently close to one, for a
sufficiently large $n$. This includes the probability to obtain a  very small $E_K$, starting from any initial error $E_{k_0}$, and to
preserve it from further perturbations. In the next, we will give the formal definition of the probability described now and we will
refer to it as the probability of  \emph{n}-\emph{step error decay}.

Before that, we need  to evaluate the detection error probability, which is defined and computed in the
next paragraph.


\subsection{Computation of the Detection Error Probability}
Let us define the stochastic process
$(D_k)_{k=0,1,\dots}$ that represents the distance between  the states estimated by the One State
procedure and the ones corresponding to the system with compensation:
  \bea  \left\{\begin{array}{l}
 D_k=\hat{X}_k-X_k=e^{\tau \MA }D_{k-1}+\frac{\hat{Z}_{k-1}-z_{k-1}}{\hat{Z}_{k-2}}~\MM_{\tau,k}\\
 D_0=0.\\
  \end{array}\right.
  \eea
Then,
\begin{definition}
Given $k\in\N$, $\dd \in \R^n$ and $\zeta\in\{\zz,\zu\}$, we define the Detection Error Probability (\emph{DEP} for short) as \bea
\text{\emph{DEP}} (k,\dd,\zeta)=P\left(\hat{Z_k}\neq z_k|D_{k}=\dd, \hat{Z}_{k-1}= \zeta\right).\eea
\end{definition}

By the definition of $D_k$, the DEP is equal to
 \be P(\hat{Z_k}\neq z_k,D_{k+1}=e^{\tau \MA}d+\frac{z_k^c-z_k}{z_{k-1}}\MM_{\tau,k+1}|D_{k}=\dd,\hat{Z}_{k-1}=z_{k-1})\ee where $z_k^c$
indicates the complementary of $z_k$ in $\{\zz,\zu\}$. This probability may be interpreted as the transition probability of the
Markov Process \bea(D_k,\hat{Z}_{k-1})_{k=0,1,\dots}\eea in the state space $\mathbf{D}\times\{\zz,\zu\}$, $\mathbf{D}\subset \R^n$,
with
starting state $(D_0,\hat{Z}_{-1})=(0,\zz)$. 

The DEP, which is fundamental to calculate the probability of the event $\{E_{k_0+n}=e^{n\tau A}E_{k_0}\}$ as shown in the next paragraph, can be analytically evaluated in the case of scalar output ($m=1$ in the system (\ref{fault})) and extended to the case $m>1$ with no particular difficulty, through some numerical techniques. In this paper, we discuss in the case $m=1$, which turns out to be interesting for the possibility of analytically describing the behavior of the DEP with respect to the parameters and to analytically derive design criteria for the fault detection. In the sequel, we then assume $Y_k,R_k \in \R$, $k=1,\dots,K$.\\

Let

\bea S_k^{w}=\MC e^{\tau \MA }\hat{X}_{k-1}+\frac{w}{\hat{Z}_{k-2}}\MC \MM_{\tau,k}\in \R\eea with $w \in\{\zz,\zu\}$ be the two
possible received signals estimated by the One State Algorithm at the generic step $k$.
The DEP is then computed in the following
\begin{proposition}\label{DEP_computation}
For any $k=1,2,\dots,K$,
\be\begin{split}
&\emph{\DEP}(k-1,\dd,\zeta)=\\
&~=\frac{1}{2}\emph{\erfc} \left(\frac{\Big{|}\frac{\zz-\zu}{2\zeta}\MC
\MM_{\tau,k}\Big{|}+\MC  e^{\tau \MA
}\dd\left[\left(1-2\mathds{1}_{\{\zz\}}(z_{k-1})\right)\left(1-2\mathds{1}_{(S_k^{\zu},+\infty)}(S_k^{\zz})\right)\right]}{\sigma\sqrt{2
}}\right)\end{split}\ee
\end{proposition}
\proof 
Under the hypothesis that $z_{k-1}=\zu$ the DEP is given by:
\bea\begin{split}
 &\DEP(k-1,\dd,\zeta)|_{(z_{k-1}=\zu)}=P\left( \hat{Z}_{k-1}=\zz\Big{|}D_{k-1}=\dd,\hat{Z}_{k-2}=
\zeta, z_{k-1}=\zu\right)\\&=P\left(|R_k-S_k^{\zz}|<|R_k-S_k^{\zu}|~~\Big{|}D_{k-1}=\dd,\hat{Z}_{k-2}=
\zeta, z_{k-1}=\zu\right)\\
&=\left\{\begin{array}{ll}
P\left(R_k<\frac{S_k^{\zu}+S_k^{\zz}}{2}~~\Big{|}D_{k-1}=\dd,\hat{Z}_{k-2}=
\zeta, z_{k-1}=\zu\right) & \text{ if } S_k^{\zu}>S_k^{\zz}\\
P\left(R_k\geq\frac{S_k^{\zu}+S_k^{\zz}}{2}~~\Big{|}D_{k-1}=\dd,\hat{Z}_{k-2}=
\zeta, z_{k-1}=\zu\right) & \text{ otherwise.}\\
         \end{array}
\right.   \end{split}\eea
If $S_k^{\zu}>S_k^{\zz}$:
\bea\begin{split}
 &P\left(R_k<\frac{S_k^{\zu}+S_k^{\zz}}{2}~~\Big{|}D_{k-1}=\dd,\hat{Z}_{k-2}=
\zeta, z_{k-1}=\zu\right)=\\&=P\left(R_k<\MC e^{\tau
\MA }\hat{X}_{k-1}+\frac{\zz+\zu}{2\zeta}\MC \MM_{\tau,k}~~|D_{k-1}=\dd\right)\\
&=P\left(\MC X_k+N_k<\MC e^{\tau
\MA }\hat{X}_{k-1}+\frac{\zz+\zu}{2\zeta}\MC \MM_{\tau,k}~~|D_{k-1}=\dd\right)\\    
&=P\left(\MC e^{\tau \MA}X_{k-1}+\frac{\zu}{\zeta}\MC \MM_{\tau,k}+N_k<\MC e^{\tau
\MA }\hat{X}_{k-1}+\frac{\zu+\zz}{2\zeta}\MC \MM_{\tau,k}~~|D_{k-1}=\dd\right)\\    
&=P\left(N_k<\MC e^{\tau
\MA }\dd+\frac{\zz-\zu}{2\zeta}\MC \MM_{\tau,k}\right)\\    
 &=\frac{1}{2}\erfc\left(\frac{-\MC  e^{\tau \MA }\dd+\frac{\zu-\zz}{2\zeta}\MC \MM_{\tau,k}}{\sigma\sqrt{2}}\right).
\end{split}\eea
The last step depends on the gaussian distribution of $N_k$; notice also that
$\frac{\zu-\zz}{\zeta}\MC \MM_{\tau,k}=S_k^{\zu}-S_k^{\zz}>0$. 

It follows also that for $S_k^{\zu}\leq S_k^{\zz}$:
\bea
 P\left(R_k\geq\frac{S_k^{\zu}+S_k^{\zz}}{2}~~\Big{|}D_{k-1}=\dd,\hat{Z}_{k-2}=
\zeta, z_{k-1}=\zu\right)=1-\frac{1}{2}\erfc\left(\frac{-\MC  e^{\tau \MA }\dd+\frac{\zu-\zz}{2\mathrm{z}}\MC
\MM_{\tau,k}}{\sigma\sqrt{2}}\right).\eea where $\frac{\zu-\zz}{\zeta}\MC
\MM_{\tau,k}=S_k^{\zu}-S_k^{\zz}\leq 0$.\\
Summing up,
\bea\begin{split}
 &\DEP(k-1,\dd,\zeta)|_{(z_{k-1}=\zu)}=\\&=P\left(|R_k-S_k^{\zz}|<|R_k-S_k^{\zu}|~~\Big{|}D_{k-1}=\dd,\hat{Z}_{k-2}=
\zeta, z_{k-1}=\zu\right)\\
&=\left\{\begin{array}{ll}
\frac{1}{2}\erfc\left(\frac{-\MC  e^{\tau \MA }\dd+\frac{\zu-\zz}{2\zeta}\MC \MM_{\tau,k}}{\sigma\sqrt{2}}\right) & \text{ if }
S_k^{\zu}>S_k^{\zz}\\
1-\frac{1}{2}\erfc\left(\frac{-\MC  e^{\tau \MA }\dd+\frac{\zu-\zz}{2\zeta}\MC \MM_{\tau,k}}{\sigma\sqrt{2}}\right)& \text{
otherwise.}\\
\end{array}\right.   \end{split}\eea

This actually corresponds to the false negative probability. The false positive probability
$\DEP(k-1,\dd,\zeta)|_{(z_{k-1}=\zz)}$  can be computed in the same way and the result is:
\bea\begin{split}
 &\DEP(k-1,\dd,\zeta)|_{(z_{k-1}=\zz)}=P\left( \hat{Z}_{k-1}=\zu\Big{|}D_{k-1}=\dd,\hat{Z}_{k-2}=
\zeta, z_{k-1}=\zz\right)\\
&=P\left(|R_k-S_k^{\zu}|<|R_k-S_k^{\zz}|~~\Big{|}D_{k-1}=\dd,\hat{Z}_{k-2}=
\zeta, z_{k-1}=\zz\right)\\
&=\left\{\begin{array}{ll}
1-\frac{1}{2}\erfc\left(\frac{-\MC  e^{\tau \MA }\dd-\frac{\zu-\zz}{2\zeta}\MC \MM_{\tau,k}}{\sigma\sqrt{2}}\right) & \text{ if }
S_k^{\zu}>S_k^{\zz}\\
\frac{1}{2}\erfc\left(\frac{-\MC  e^{\tau \MA }\dd-\frac{\zu-\zz}{2\zeta}\MC \MM_{\tau,k}}{\sigma\sqrt{2}}\right)& \text{
otherwise.}\\
\end{array}\right.   \end{split}\eea
The thesis is then proved.\qed
\begin{remark}
If $\dd=0\in \R^n$,
\be\label{remk}\begin{split}
\emph{\DEP}(k-1,0,\zeta)&=\frac{1}{2}\emph{\erfc}
\left(\frac{\Big{|}\frac{\zz-\zu}{2\zeta}\MC
\MM_{\tau,k}\Big{|}}{\sigma\sqrt{2}}\right)\\&=
\frac{1}{2}\text{\emph{erfc}}\left(\frac{|S_k^{\zz}-S_k^{\zu}|/2}{\sigma\sqrt{2}}\right).\end{split}\ee
This expression suggests an Information theoretic intepretation of our problem. In fact, the presence of the gaussian noise in the data
lecture
can be thought as if signal $y_k$  was transmitted on an Additive White Gaussian Noise (AWGN) channel. If $D_{k-1}=0$, $y_k$ can be
 $S_k^{\zz}$ or $S_k^{\zu}$.  Moreover, if we shift the
signals by their average, so that they become antipodal $\pm \frac{S_k^{\zz}-S_k^{\zu}}{2}$, the average energy per channel use at step
$k$ is $\mathcal{E}_k=\left(\frac{S_k^{\zz}-S_k^{\zu}}{2}\right)^2$. Given that the spectral density of the gaussian noise is
$N_0=2\sigma^2$, the argument of the erfc function in (\ref{remk}) turns out to be the square root of the so called Signal-to-Noise
Ratio (SNR), defined as
$\SNR_k=\mathcal{E}_k/N_0$, of our ideal channel.

Generally, the SNR compares the magnitudes of the transmitted signal and of the channel noise and 
it is widely used in Informatiom Theory to describe channel performance. In our framework, the SNR determines the 
reliability of the detection, say the reliability of the  channel where $y_k$ is ideally transmitted.
This remark emphasizes that our problem is analogous to a common digital-transmission paradigm and bears out the idea of using
decoding techniques to the detection task.

In the next, we will use the common dB notation for the SNR, that is, we express it as $10\log_{10}$ of its value.
\end{remark}
\begin{remark}
Since typically $\zu<\zz$, by expression (\ref{remk}) we have
\bea 
\emph{\DEP}(k-1,0,\zu)<\emph{\DEP}(k-1,0,\zz)
.\eea
Given that $\hat{Z}_{k-2}=\zu$ is generally more likely when $z_{k-2}=\zu$ (otherwise our detection method would be improper),
we can conclude that our detection algorithm is more reliable after the failure, or, in other terms, it is more sensitive to false
positives.
\end{remark}

\subsection{Computation of the Probability of \emph{n}-step Error Decay}
Given a time interval $[k_0,k_0+n)$, $k_0,n \in \N$, $k_0\geq 1$, we can formally define  the probability of
\emph{n}-step error decay (EDP$^n$ for short) as
\bea
\begin{split}
 & \EDP^n (k_0,\dd,\zeta,\eta)=\\
&P\left(E_{k_0+n}=e^{n\tau
\MA}E_{k_0}\big{|}D_{k_0-1}=\dd, \hat{Z}_{k_0-2}=\zeta, z_k=\eta\text{ for any } k=k_0-1,\dots,k_0+n-1)\right)\\
\end{split}\eea
where $\dd\in \R^n$, $\zeta,\eta \in\{\zz,\zu\}$. Notice that $z_k$ is assumed to be constant in $[k_0-1,k_0+n-1]$, that is, we consider
the system before or after a failure event. Recalling the Proposition \ref{equiv_events}, the EDP is connected to the DEP by the following expression:
\bea
\begin{split}
 & \EDP^1 (k_0,\dd,\zeta,\eta)=P\left(E_{k_0+1}=e^{\tau
\MA}E_{k_0}\big{|}D_{k_0-1}=\dd, \hat{Z}_{k_0-2}=\zeta, z_{k_0-1}=z_{k_0}=\eta\right)\\
&=P\left(\hat{Z}_{k_0-1}=z_{k_0}\big{|}D_{k_0-1}=\dd, \hat{Z}_{k_0-2}=\zeta, z_{k_0-1}=z_{k_0}=\eta\right)\\
&=1-\DEP(k_0-1,\dd,\zeta)_{\big{|}z_{k_0-1=\eta}}\\
\end{split}\eea
that is, the Error decays when the detection is correct. Notice that this relation between EDP and DEP subsists in virtue of the condition $z_{k_0-1}=z_{k_0}$: if $k_0$ were a switch point, the feedback delay would produce a deviation in the Error Function in case of correct detection.\\
Generalizing to $n$ steps,
\bea\begin{split}& \EDP^n(k_0,\dd,\zeta,\eta)=\\&=P(\hat{Z}_{k_0-1}=\hat{Z}_{k_0}=\dots=\hat{Z}_{k_0+n-2}
=\eta\big{|}D_ {k_0-1 } =\mathrm{d},\hat{Z}_{k_0-2}=\zeta)\\
&=P\left( (D_{k_0},\hat{Z}_{k_0-1})=(e^{\tau
\MA}\mathrm{d},\eta)|(D_{k_0-1},\hat{Z}_{k_0-2})=(\mathrm{d},\zeta)\right)\cdot\\
&~~\cdot\prod_{m=1}^{n-1}P\left( (D_{k_0+m},\hat{Z}_{k_0+m-1})=(e^{(m+1)\tau
\MA}\mathrm{d},\eta)\big{|}(D_{k_0+m-1},\hat{Z}_{k_0+m-2})=(e^{m\tau
\MA}\mathrm{d},\eta)\right)\\
&=\EDP^1(k_0,\dd,\zeta,\eta)\prod_{m=1}^{n-1}\EDP^{1}(k_0+m,e^{m\tau\MA}\dd,\eta,\eta)\\
&=\big{(}1-\DEP(k_0-1,\dd,\zeta)\big{)}_{\big{|}z_{k_0-1}=\eta}\prod_{m=1}^{n-1}\big{(}1-\DEP(k_0+m-1,e^{m\tau
A}\dd,\eta)\big{)}_{\big{|}z_{k_0+m-1}=\eta}\\
\end{split}\eea
By  Proposition \ref{DEP_computation}, this is equal to
\be\label{nsteps_byprop}\begin{split}
& \EDP^n(k_0,\dd,\zeta,\eta)=\\
&=\frac{1}{2}\erfc \left(-\frac{\Big{|}\frac{\zz-\zu}{2\zeta}\MC
\MM_{\tau,k_0}\Big{|}+\MC  e^{\tau \MA
}\dd\left[\left(1-2\mathds{1}_{\{\zz\}}(\eta)\right)\left(1-2\mathds{1}_{(S_k^{\zu},+\infty)}(S_k^{\zz})\right)\right]}{\sigma\sqrt
{ 2
}}\right)\\
&\cdot\prod_{m=1}^{n-1}\frac{1}{2}\erfc \left(-\frac{\Big{|}\frac{\zz-\zu}{2\eta}\MC
\MM_{\tau,k_0+m}\Big{|}+\MC  e^{(m+1)\tau \MA
}\dd\left[\left(1-2\mathds{1}_{\{\zz\}}(\eta)\right)\left(1-2\mathds{1}_{(S_{k+m}^{\zu},+\infty)}(S_{k+m}^{\zz})\right)\right]}{
\sigma\sqrt{ 2}}\right).\end{split}\ee
Our next goal is to evaluate the $\EDP^n$ in different instances of system (\ref{fault},\ref{failure}). First of all, let us
distinguish what happens before and after the failure.
\subsection{False positive evaluation}
Let suppose the system to be affected by a failure according to the model (\ref{failure}) with $k_F\geq 1$, that is, the system is not
faulty from the beginning. In particular, since there is no compensation at the first time step (or equivalently $\hat{Z}_{-1}=\zz$), no
false positive is produced at $k=0$. Then, studying the $\EDP$ in $[1,k_F)$ actually
corresponds to evaluate the probability that no false postives occur during the whole pre-failure transient regime. Given that $D_0=0$,
we have
\be\label{false_positive}\begin{split}
& \EDP^{k_F-1}(1,0,\zz,\zz)=\prod_{m=1}^{k_F-1}\frac{1}{2}\erfc \left(-\frac{\Big{|}\frac{\zz-\zu}{2\zz}\MC
\MM_{\tau,m}\Big{|}}{\sigma\sqrt{ 2}}\right).\end{split}\ee
Since $E_1=0$ and $D_0=0$, then $\EDP^{k_F-1}(1,0,\zz,\zz)=P(E_{k_F}=0)=P(D_{k_F}=0)$.
\subsection{Switch Point}
Suppose that $D_{k_F}=0$, then in particular, $\hat{Z}_{k_F-1}=z_{k_F-1}$ and $\hat{Z}_{k_F-1}\neq z_{k_F}$. In other terms, the detection is correct, but the compensation, based on the detection at the previous step, is not efficient in correspondance of a switch point. Our detection method cannot control what happens at step at step $k_F$, that is, in the time interval $[T_F,T_F+\tau)$.

\subsection{False negative evaluation}
Given that we cannot control the system immediately after the switch point, it is likely that $E_{k_F+1}\neq 0$.  We now want to study the probability of decay of the Error Function towards zero, which actually corresponds to the evaluation of the false negatives. In fact, under the hypothesis $D_{k_F}=0$ (i.e., no false positives and in particular $\hat{Z}_{k_F-1}=\zz$), for any $n\in \N$,
\be\label{nsteps_byprop_d0}
\begin{split}
&\EDP^n(k_F+1,0,\zz,\zu)=\EDP^1(k_F+1,0,\zz,\zu)\prod_{m=1}^{n-1}\EDP^1(k_F+1+m,0,\zu,\zu)\\
&\frac{1}{2}\erfc \left(-\frac{\Big{|}\frac{\zz-\zu}{2\zu}\MC
\MM_{\tau,k_F+1}\Big{|}}{\sigma\sqrt{2}}\right)
\prod_{m=1}^{n-1}\frac{1}{2}\erfc \left(-\frac{\Big{|}\frac{\zz-\zu}{2\zu}\MC
\MM_{\tau,k_F+1+m}\Big{|}}{\sigma\sqrt{2}}\right).  
\end{split}
\ee
Notice that $n$ can be any positive integer, since the failure state is not reversible. Moreover, it is clear that if $n\to\infty$, then $EDP^n\to 0$, that is, it is not likely that the Error decays to zero and remains null forever. However, we can approximate this ideal situation, as we will see in the next.

The considerations about the EDP made in this section are now applied to the case of constant input $f(t)$. More precisely we will exploit them to establish suitable design criteria, that is, which is the best choice of parameters to obtain the maximum performance from the One State Algorithm.

\subsection{Constant input $f(t)$ }\label{constf}
If the input $f(t)$ is constant, say $f\equiv 1$, the system evolution does not depend on time step $k$. In fact,
$\MM_{\tau,k}=\MM_{\tau}=
(e^{\tau \MA }-\mathds{I})\MA ^{-1}\MB$ for any $k=1,\dots,K$.
Hence, 
\be\begin{split} &\EDP^{n}(1,0,\zz,\zz)=\left[\frac{1}{2}\erfc \left(-\frac{\Big{|}\frac{\zz-\zu}{2\zz}\MC
\MM_{\tau}\Big{|}}{\sigma\sqrt{2}}\right)\right]^{n}\end{split}\ee
for any $n \in \N$ such that $n+1\leq k_F$ and
\be\begin{split} &\EDP^{n}(k_F+1,0,\zz,\zu)=\frac{1}{2}\erfc \left(-\frac{\Big{|}\frac{\zz-\zu}{2\zz}\MC
\MM_{\tau}\Big{|}}{\sigma\sqrt{2}}\right)\left[\frac{1}{2}\erfc \left(-\frac{\Big{|}\frac{\zz-\zu}{2\zu}\MC
\MM_{\tau}\Big{|}}{\sigma\sqrt{2}}\right)\right]^{n-1}.\end{split}\ee
In terms of signal-to-noise ratio, we can write
 \bea \sqrt{\SNR(\eta)}=\frac{\left |
\frac{\zu-\zz}{2\eta}\MC \MM_{\tau}\right |}{\sigma\sqrt{2}} \eea 
so that
\bea
\begin{split}
&\EDP^{n}(1,0,\zz,\zz)=\left[\frac{1}{2}\erfc \left(-\sqrt{\SNR(\zz)}\right)\right]^{n}  \\
&\EDP^n(k_F+1,0,\zz,\zu)=\frac{1}{2}\erfc \left(-\sqrt{\SNR(\zz)}\right)\left[\frac{1}{2}\erfc \left(\sqrt{\SNR(\zu)}\right)\right]^{n-1}.\\
\end{split}\eea
Under the hypothesis $0<\zu<\zz=1$, $\SNR(\zz)<\SNR(\zu)$, that is $\EDP^m(k_0,0,\zz,\zz)<\EDP^m(k_1,0,\zu,\zu)$; in other terms, our
detection algorithm is more sensitive to false positives, then our fault tolerant control method is more efficient \emph{after} the
failure. Hence, the suitable design criteria for the pre-failure state will automatically be appropriate also for the post-failure
state. This is why in the next we will generically name
 \be\label{simply}\SNR=\SNR(\zz)~~\text{ and }~~ \EDP^n=\EDP^n(k_0,0,\zz,\zz)=\left[\frac{1}{2}\erfc
\left(-\sqrt{\SNR}\right)\right]^{n}.\ee

The next section is devoted to the study of design criteria for our FTC system, on the basis of the theoretical analysis
developed in the last pages. Particular attention will be paid to the case of constant $f(t)$, for which optimal criteria can be
formulated.
\section{Design Criteria}\label{descri}
In this section, our aim is to provide the design criteria to obtain the best performance from our FTC scheme, based on
the One State Algorithm. 

The key point of this issue is that the controller is supposed to be free to choose the sampling time step $\tau$, hence our goal is to
give the criteria to determine the \emph{otpimal} $\tau$, which, in our framework, can be defined as the one  that
\emph{minimizes} the Error Function, in the sense that we now explain. Given the failure system (\ref{fault},\ref{failure}) and a time
window $W=n\tau$ not containing the switch point, our first purpose is to maximize the probability that $E_k$ remains null (if we set
before the failure) or decays to zero (if we set after the failure) along the interval $W$. Furthermore, given that in $(T_F,T_F+\tau]$
a correct detection causes a failed compensation and a consequent abrupt deviation in the output $y$ (as we will show in the numerical
simulations), our second purpose is to minimize the peak of this unavoidable deviation.

This qualitative discussion is now quantified in two different input instances: $f(t)$ constant and $f(t)$ sinusoidal. As far as the
first case in concerned, we will show that the theoretic analysis of Section \ref{theoanalysis} provides the instrument to determine the
sampling time that minimizes the Error Function in an analytic way. On the other hand, when the input is not constant some
difficulties arise in the definition of the optimal $\tau$; however, we will explain how to obtain suitable values of $\tau$ by a
numerical numerical computation, still based on the analysis of Section \ref{theoanalysis}.

\subsection{Design Criteria in the case of constant input $f(t)$}
Recalling the Paragraph \ref{constf} and in particular the simplified notation (\ref{simply}), let us explain how to define the optimal
$\tau$ when $f(t)\equiv 1$. As just said, we aim to maximize the EDP in a given time window $W$ not containing the failure instant and
to minimize the peak of the deviation immediately after the failure. In particular, if $E_{k_F}=0$, by definition \ref{defErr}, the
extent of the peak in the output is given by $\max_{t\in (0,\tau]}|\frac{\zeta_1-\zeta_0}{\zeta_0}\MC\MM_{t}|$. In brief, we intend to
provide
\be
\begin{split}
&\tau_1=\underset{\tau>0}{\operatorname{argmax}}~\EDP^{W/\tau}~~~\text{and}~~~\tau_2=\underset{\tau>0}{\operatorname{argmin}}
~\left(\max_ { t \in (0,\tau]}|\MC\MM_{t}|\right)
\end{split}\ee
The optimum will be $\tau_1=\tau_2$, but in general this is not the case. Then, we define the optimal $\tau$ as follows: we do not look
for the maximum EDP, but we just require $\EDP^{W/\tau}>1-\varepsilon$ where $\varepsilon<<1$ is a fixed tolerance. In other
terms, we demand that the EDP be very close to 1. Then, the optimal $\tau$, indicated by
$\tau_{\text{opt}}=\tau_{\text{opt}}(\varepsilon)$, is :
\be\label{tauopt} \tau_{\text{opt}}=\underset{\tau:~\EDP^{W/\tau}>1-\varepsilon}{\operatorname{argmin}}\left(\max_{t\in(0,\tau]}|\MC\MM_t|\right). \ee
\subsubsection{Application  to the Flight Control Problem}\label{constflight}
Let us now compute $\tau_{\text{opt}}$ for the Flight Control Problem introduced in the Paragraph \ref{example}, in the case of
constant input $f(t)$.
\begin{figure}
\centering
\includegraphics[width=12cm,viewport=150 120 750 500]{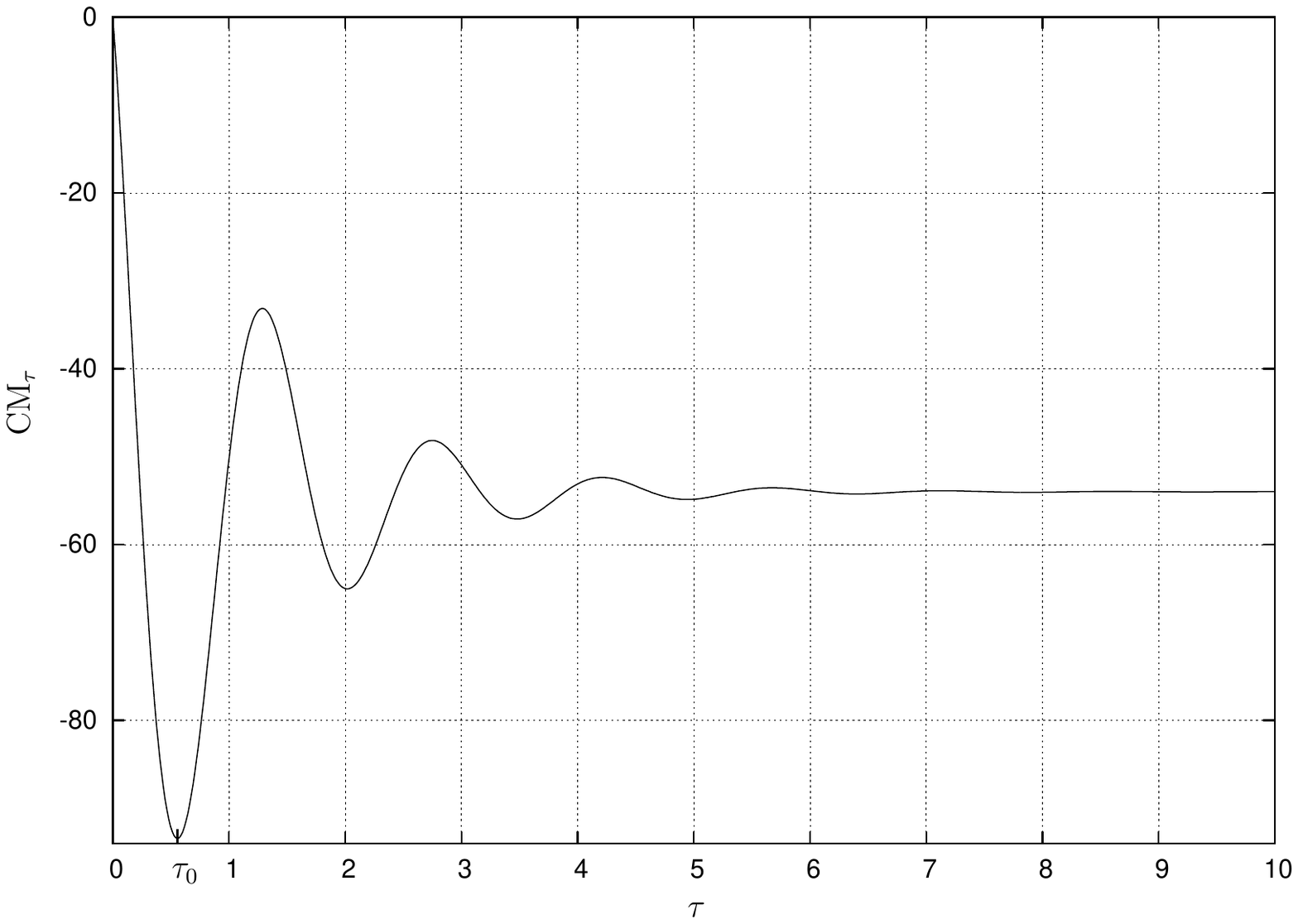}
\caption{$m_{\tau}$}\label{mtau}
\end{figure}
In the Figure \ref{mtau}, the graph of $\MC\MM_{\tau}$ in function of $\tau$ is shown. In particular, we notice that  $\MC\MM_{\tau}$ is negative for any $\tau>0$, achieves a global  minimun at $\tau_0=0.55$ and converges to a constant value for a sufficienlty large $\tau$. Then, if  $\tau>\tau_0$, $\max_{t\in(0,\tau]}|\MC\MM_t|=|\MC\MM_{\tau_0}|$, that is, the peak is fixed and we cannot control it. This undesired occurrence can be prevented by imposing \bea \tau\in (0,\tau_0]. \eea
In this interval, $\MC\MM_{\tau}$ is monotone decreasing and  $\max_{t\in(0,\tau]}|\MC\MM_t|=|\MC\MM_\tau|$.\\
Then, fixed the tolerance $\varepsilon$, our aim is the computation of
\be \tau_{\text{opt}}=\underset{\tau\in(0,\tau_0]:\EDP^{W/\tau}>1-\varepsilon}{\operatorname{argmin}}|\MC\MM_{\tau}|. \ee

Notice that \bea\EDP^{W/\tau}=\left[\frac{1}{2}\erfc \left(-\sqrt{\SNR}\right)\right]^{W/\tau}=\left[\frac{1}{2}\erfc \left(-\frac{|\frac{\zeta_1-\zeta_0}{2\zeta_0}\MC\MM_{\tau}|}{\sigma\sqrt{2}}\right)\right]^{W/\tau}\eea
is monotone increasing as a function of $\tau$. Then, let $\tau_m=\tau_m(\varepsilon)$ be the minimum $\tau$ in $(0,\tau_0]$ such that $\EDP^{W/\tau}>1-\varepsilon$ (if it exists). Then
\be \tau_{\text{opt}}=\underset{\tau\geq\tau_m}{\operatorname{argmin}}|\MC\MM_{\tau}|=\tau_m. \ee
Now, let assign numerical values to the parameter and solve the corresponding instance. Suppose that:
\be\label{numbers}
\begin{split}
&\zz=1~~~\zu=\frac{1}{2}~~~\sigma^2=2\\
&\varepsilon=10^-3~~~W=20\\
\end{split}\ee
\begin{figure}
\centering
\includegraphics[width=10.5cm,viewport=150 120 750 500 ]{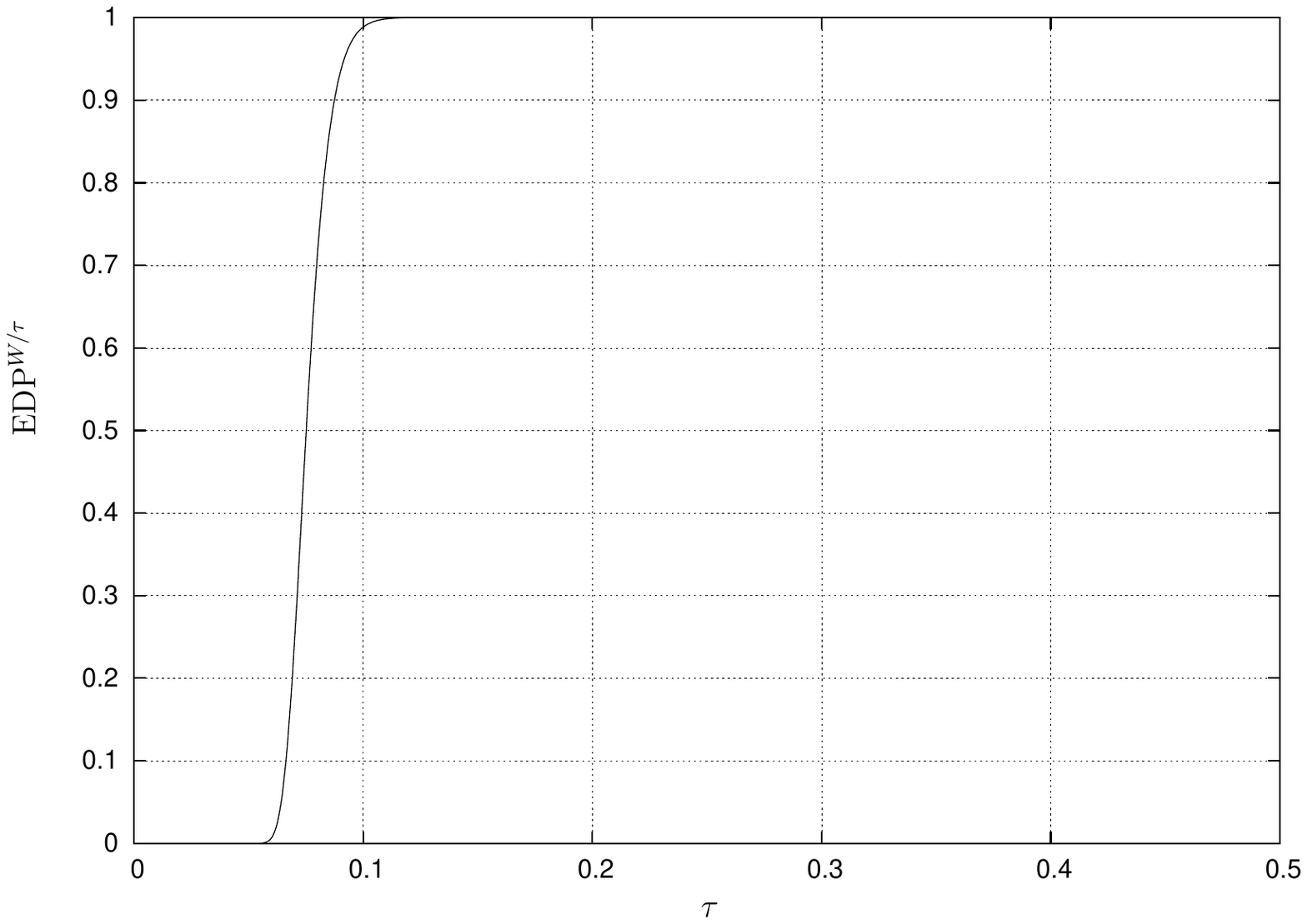}
\includegraphics[width=10.5cm,viewport=150 120 750 500 ]{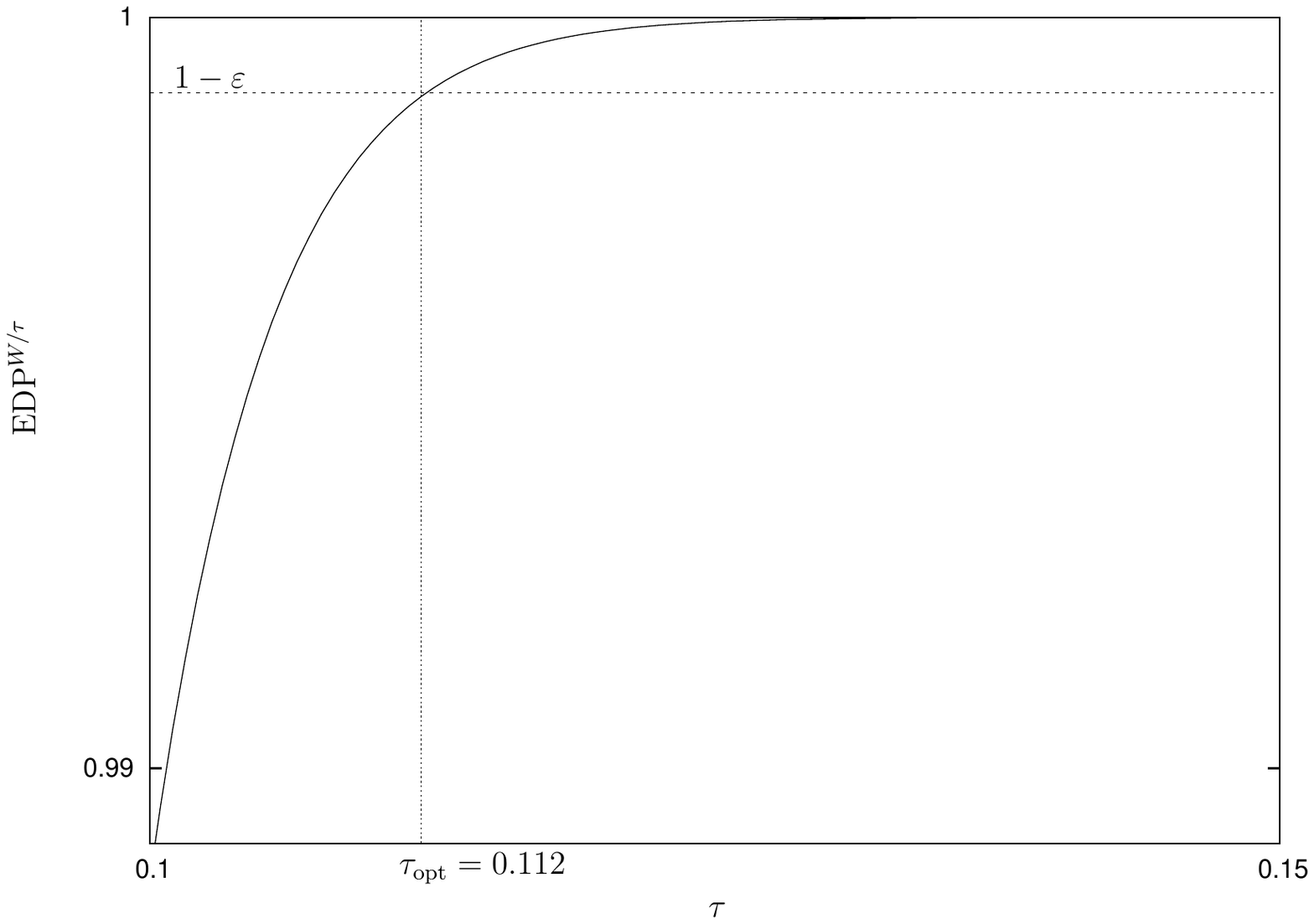}
\caption{$EDP^{W/\tau}$ in function of $\tau$ in the instance (\ref{numbers}). The second graph is a zoom that allows to see that $\tau_{\text{opt}}=0.112$}\label{edp}
\end{figure}
In this case, $\tau_{\text{opt}}=0.112$ as shown in Figure \ref{edp}.\\
\begin{figure}
\centering
\includegraphics[width=12cm,viewport=150 120 750 500]{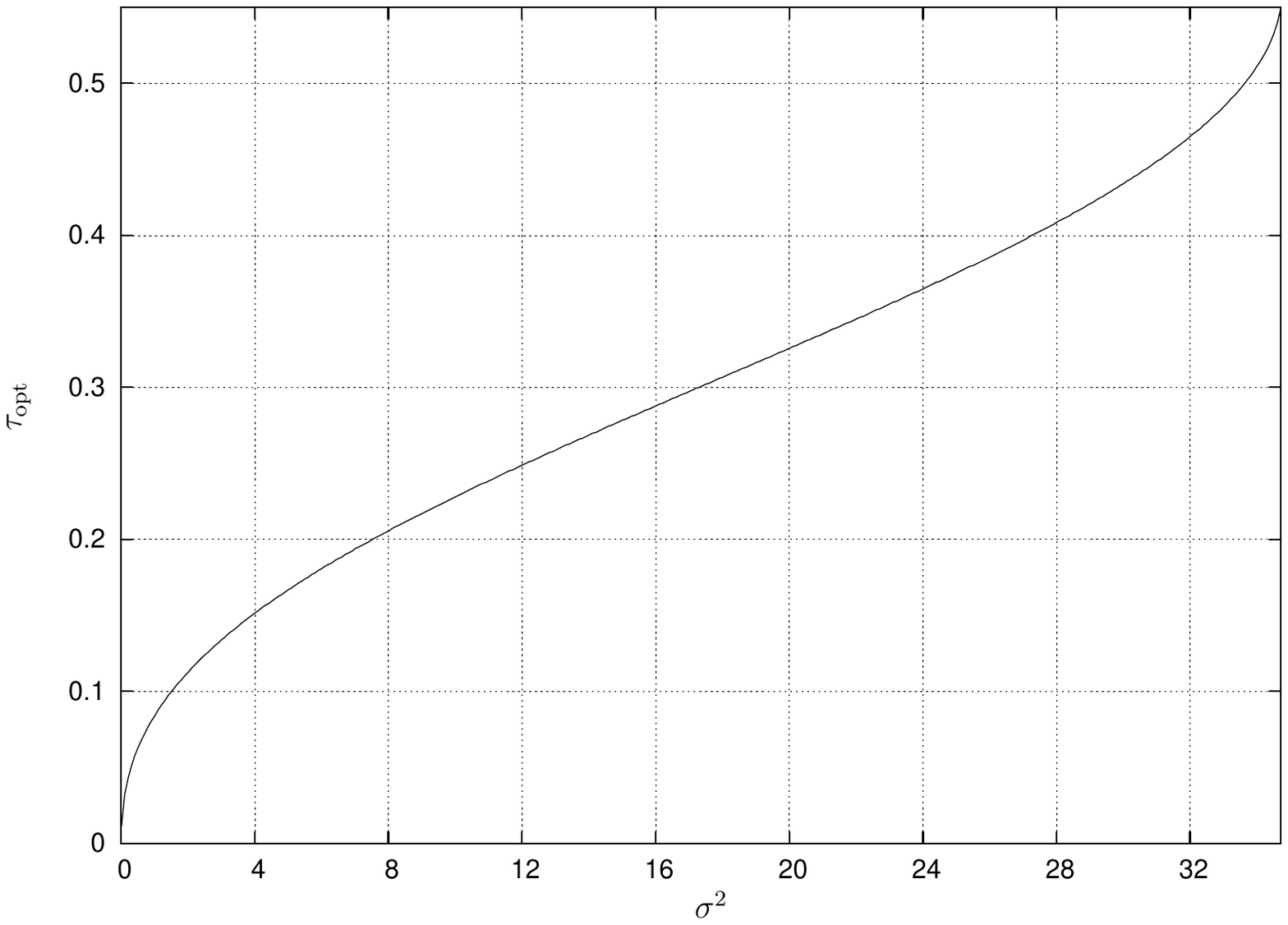}
\caption{The optimal $\tau$'s as the noise variance $\sigma^2$ changes ($\zz=1, \zu=\frac{1}{2}, \varepsilon=10^-3, W=20$)}\label{noisetau}
\end{figure}
The value of $\tau_{\text{opt}}$ clearly depends on the noise and in particular  there can exist noise values for which there is no $\tau$ making $EDP^{W/\tau}>1-\varepsilon$: for instance, this occurs if we consider $\sigma^2>34.72$ in the example (\ref{numbers}) (the range of admittible $\sigma^2$'s with the corresponding $\tau_{\text{opt}}$'s is shown in Figure \ref{noisetau}). In such situation, one should allow a lower threshold $1-\varepsilon$.\\

In Section \ref{SecSim} we will show a few simulations about the Flight Example.

\subsection{Design Criteria in the case of  input $f(t)=\sin t$}
When $f(t)$ is not constant, it is more difficult to study analytical design criteria as the quality of the detection depends on time.
In particular, at each time step $k\tau$ the detection is affected by the values of $f(t)$, $t \in ((k-1)\tau,k\tau)$, then any
detection step is different from the others and an analogous of (\ref{tauopt}) cannot be provided: roughly speaking, the optimum would
be to change $\tau$ according to the shape of $f(t)$ in each considered interval.\\
When $f(t)$ is periodic, we can suggest some numerical computation in order to fix a suitable $\tau$. In fact, if we compute
$\EDP^{W/\tau}(1,0,\zz,\zz)$ for a sufficiently large $W$, we get an idea about the sampling times that are more suitable. On the other
hand, there is no way to control the amplitude of the deviation in case of failure, since this again depends on time.
The idea is then to choose as samling time that maximises $\EDP^{W/\tau}(1,0,\zz,\zz)$ or that makes it larger than a given threshold,
being conscious that this does not arrange the anavoidable deviation.
\begin{figure}
\centering
\includegraphics[width=12cm,viewport=150 120 750 500 ]{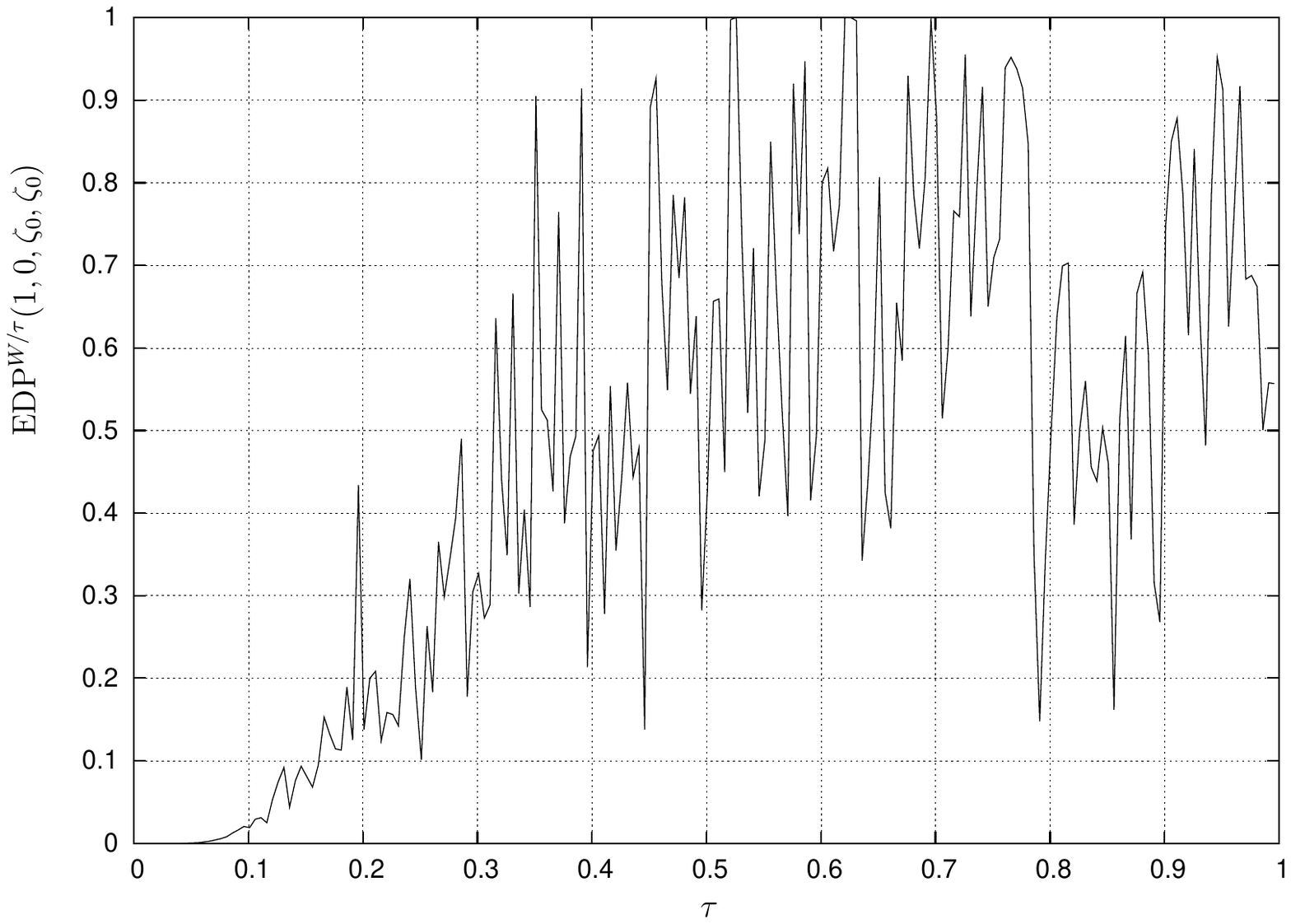}
\caption{$EDP^{W/\tau}(1,0,\zz,\zz)$ in function of $\tau$ in the instance (\ref{numbers}) ($\zz=1,\zu=\frac{1}{2}, \sigma^2=2, W=20$).}\label{edpsin}
\end{figure} 
Let us  illustrate these observations in the Flight Control Problem with $f(t)=\sin t$ and parameters given by (\ref{numbers}).
First, let us numerically compute $\EDP^{W/\tau}(1,0,\zz,\zz)$  in
function of $\tau$,  the result being presented in Figure \ref{edpsin}: the graph shows a clear unsettled behavior which cannot be
described analytically. However, it also suggests the
values of $\tau$ that give an high $\EDP^{W/\tau}(1,0,\zz,\zz)$ and which can  then considered suitable.
No general consideration can be derived, except that a very small $\tau$ is in general  not
preferable.

More details about this instance can be retrieved in the simulations presented in the next Section.

\section{Flight Control Problem: a few simulations}\label{SecSim}
In this section, we show some simulations concerning the application of the One State Algorithm to the Flight FTC example presented in the Paragraph \ref{example} and studied in the previous paragraphs.\\

In a time interval $[0,T]=[0,40]$, we suppose that a failure occurs at $T_F=20$ and causes the switch of the disturbance function
$z(t)$ from $\zz=1$ to $\zu=1/2$ ($\zu=1/2$ might represent a loss of effectiveness of $50\%$ of the elevator of the aircraft). The
lecture noise is  a gaussian random variable $\mathcal{N}(0,2)$. We consider  boht
the cases of input $f\equiv 1$ and $f(t)=\sin t$ and we show the behavior of the One State procedure for different values of $\tau$.
The graphs represent the output $y(t)$ of the system.

\begin{figure}[ht]
\centering
\includegraphics[width=6cm,viewport=225 110 680 450 ]{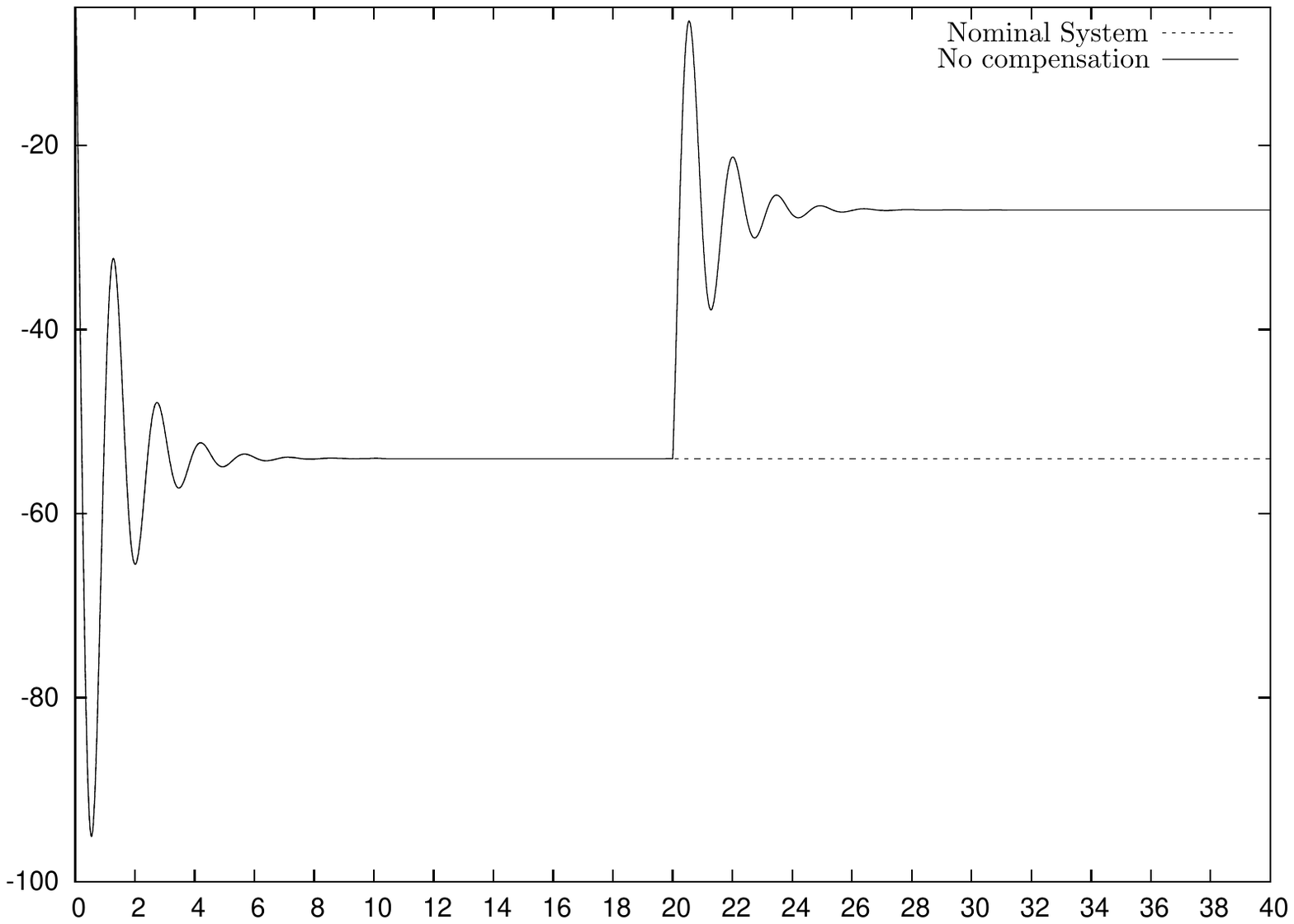}
\includegraphics[width=6cm,viewport=225 110 680 450]{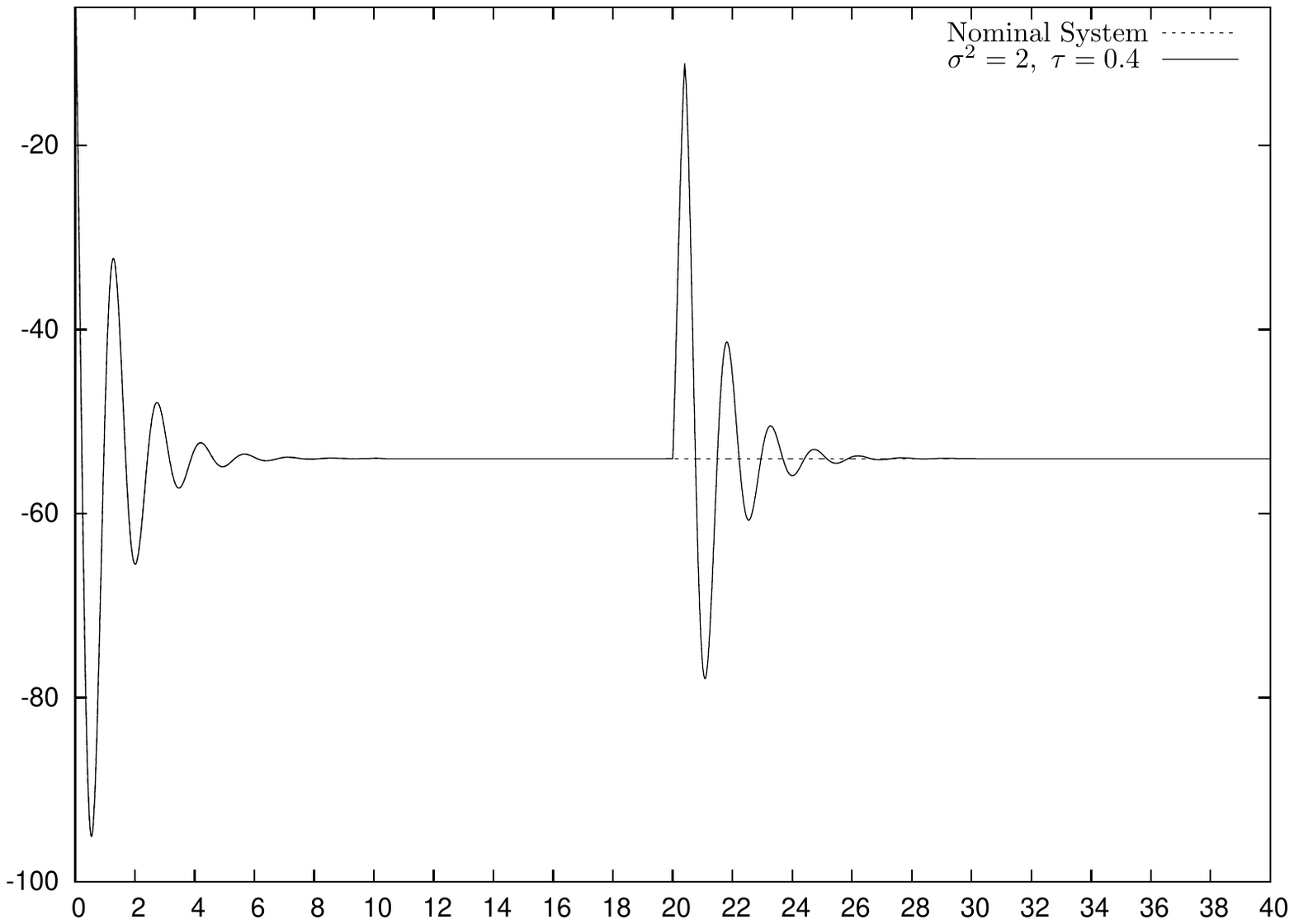}
\includegraphics[width=6cm,viewport=225 110 680 450]{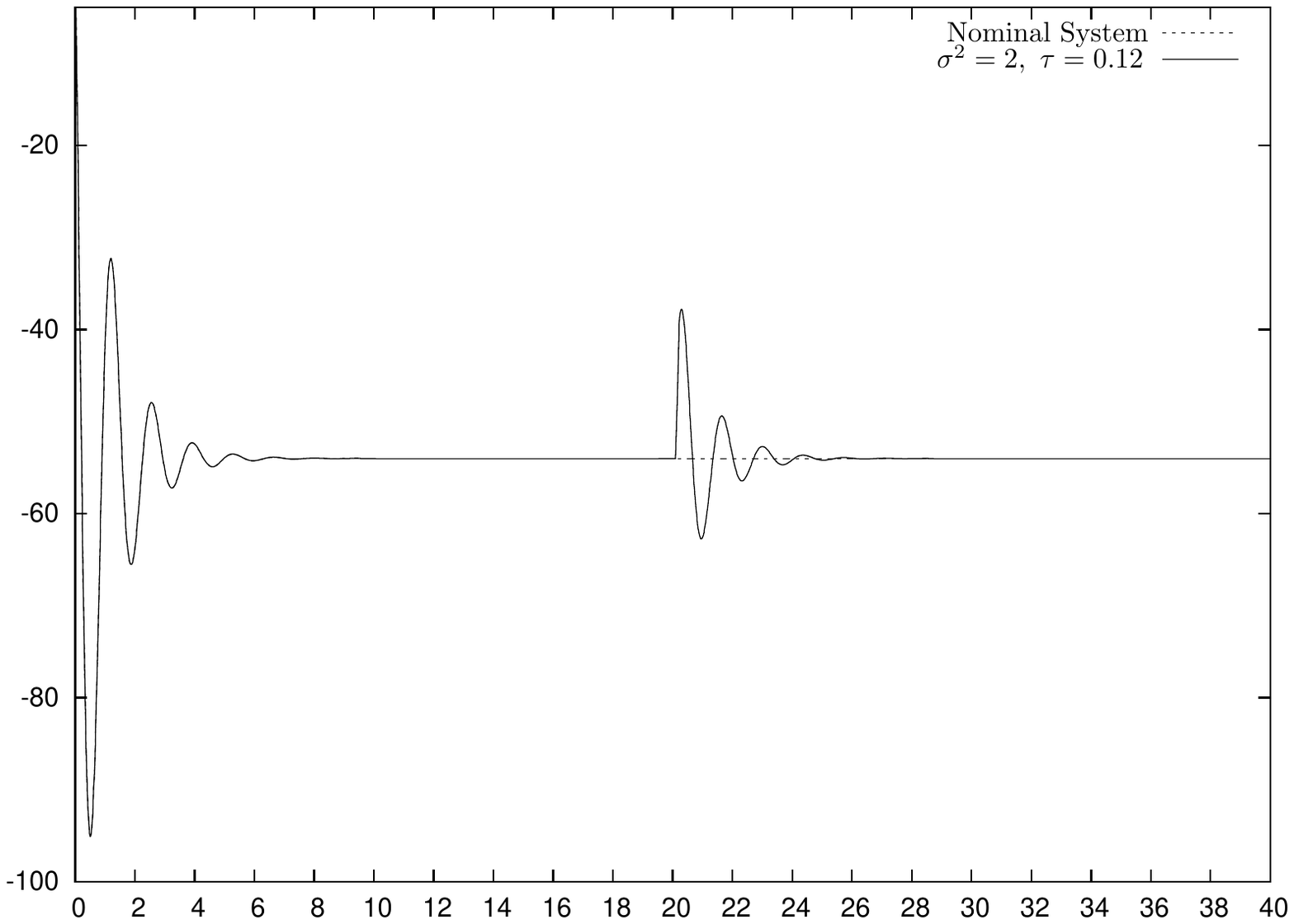}
\includegraphics[width=6cm,viewport=225 110 680 450]{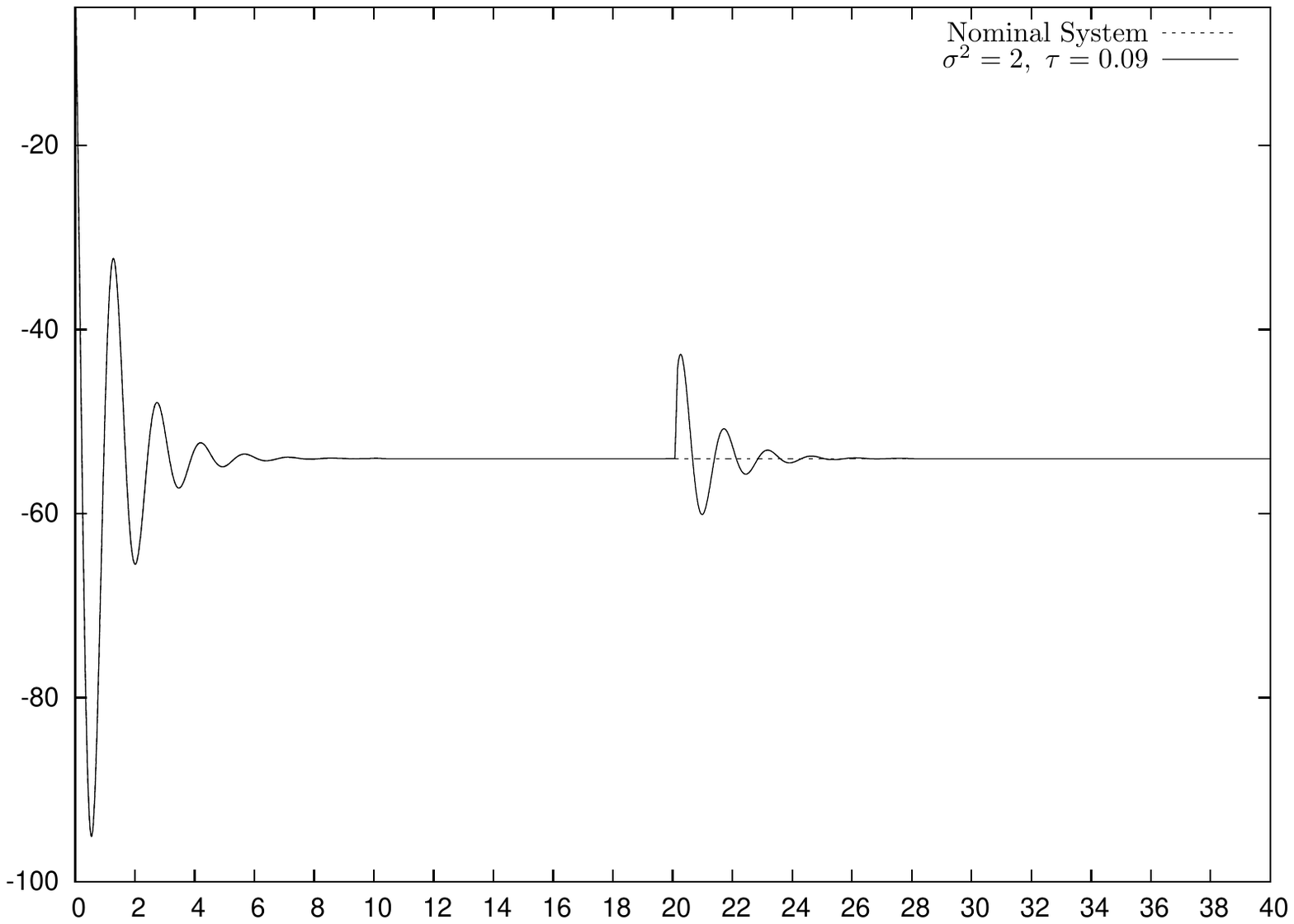}
\includegraphics[width=6cm,viewport=225 110 680 450]{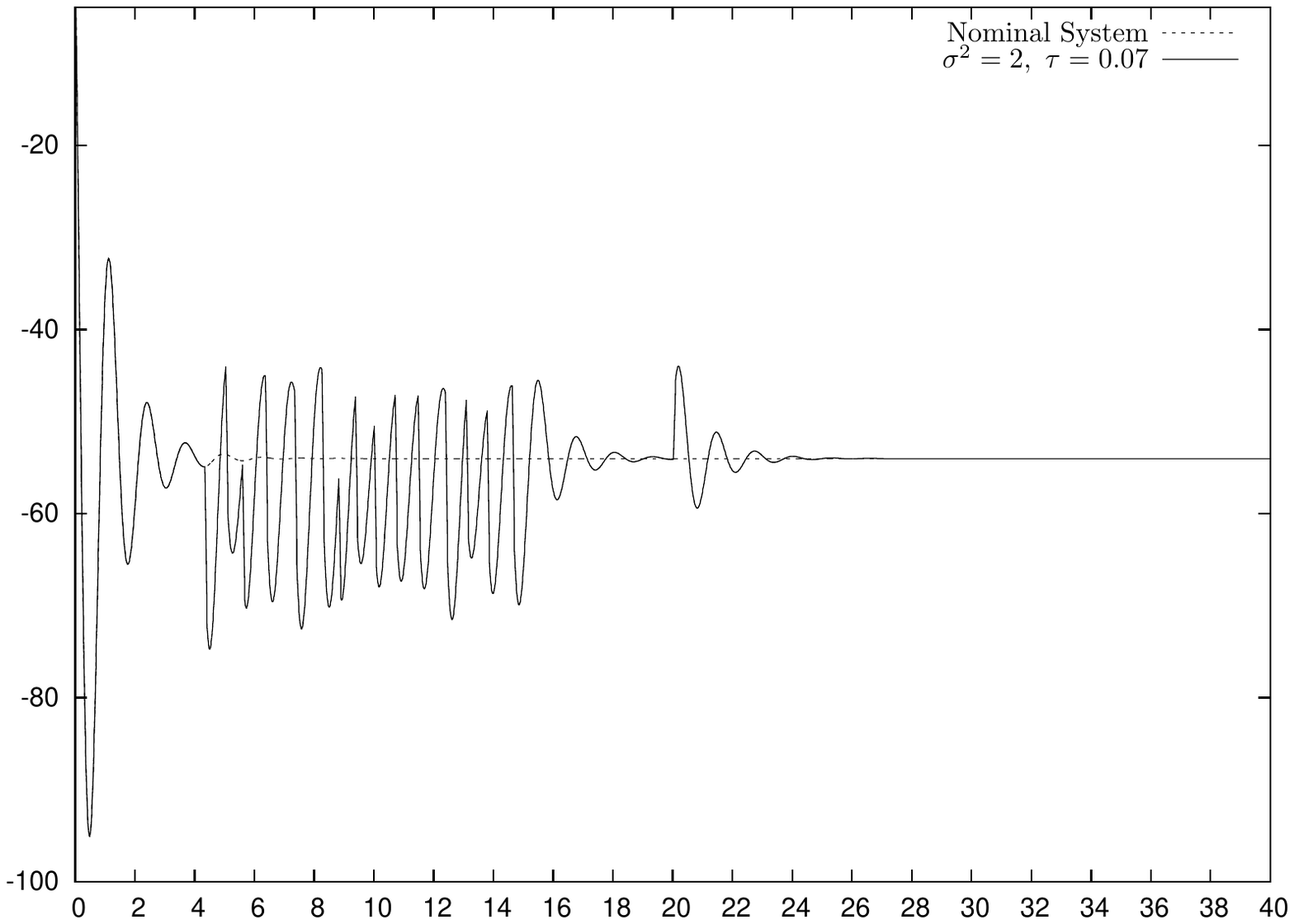}
\includegraphics[width=6cm,viewport=225 110 680 450]{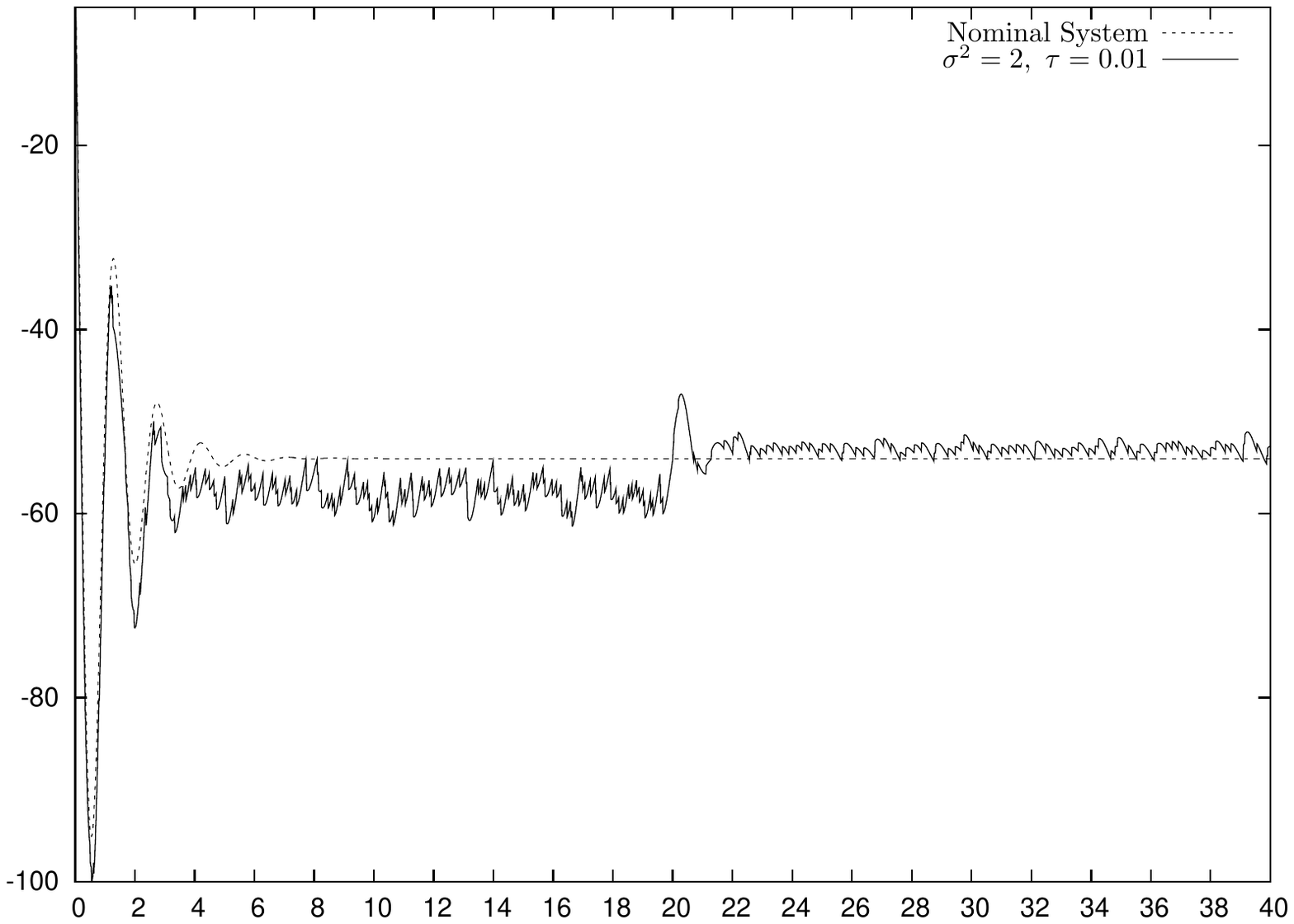}
\caption{Output y(t): Nominal System vs System with a failure at $T_F=20$, with lecture noise of variance $\sigma^2=2$ and $f\equiv 1$.
Six different cases are shown: the first graph represents the system with no control and compensation; the other ones are with
compensation, respectively with time step $\tau$ equal to $0.4$, $0.12$, $0.09$, $0.07$, $0.01$}\label{F1}
\end{figure}
Figure \ref{F1} reproduces the case $f\equiv 1$. The first graph compares the nominal system, that is, the desirable trajectory,
to the faulty system with no compensation: after the failure, the trajectory of the latter is sensibily uncorrect. In the other graphs,
we introduce the compensation using the One State Algorithm: as proved in the Paragraph \ref{constflight} ,  $\tau_{\text{opt}}=0.112$.
In the second graph, we fix $\tau=0.4$, which is larger than $\tau_{\text{opt}}$: we obtain a correct detection at each step, but the
unavoidable deviation is not optimized: in fact, considering $\tau_{\text{opt}}$ (third graph), we have a smaller peak after the
failure. Furthermore, we see that also $\tau=0.09$ is suitable, even if, the corresponding $\EDP^{W/\tau}>1-\varepsilon$. On the other
hand, $\tau=0.07$ assures a good detection only after the failure (this is consistent with our observation about the different
sensitivity ot false positives and false negatives), while a too small sampling time ($\tau=0.001)$ causes instabililty: the detection
is not reliable and the Error is always nonnull.

\begin{figure}[h]
\centering
\includegraphics[width=6cm,viewport=225 110 680 450 ]{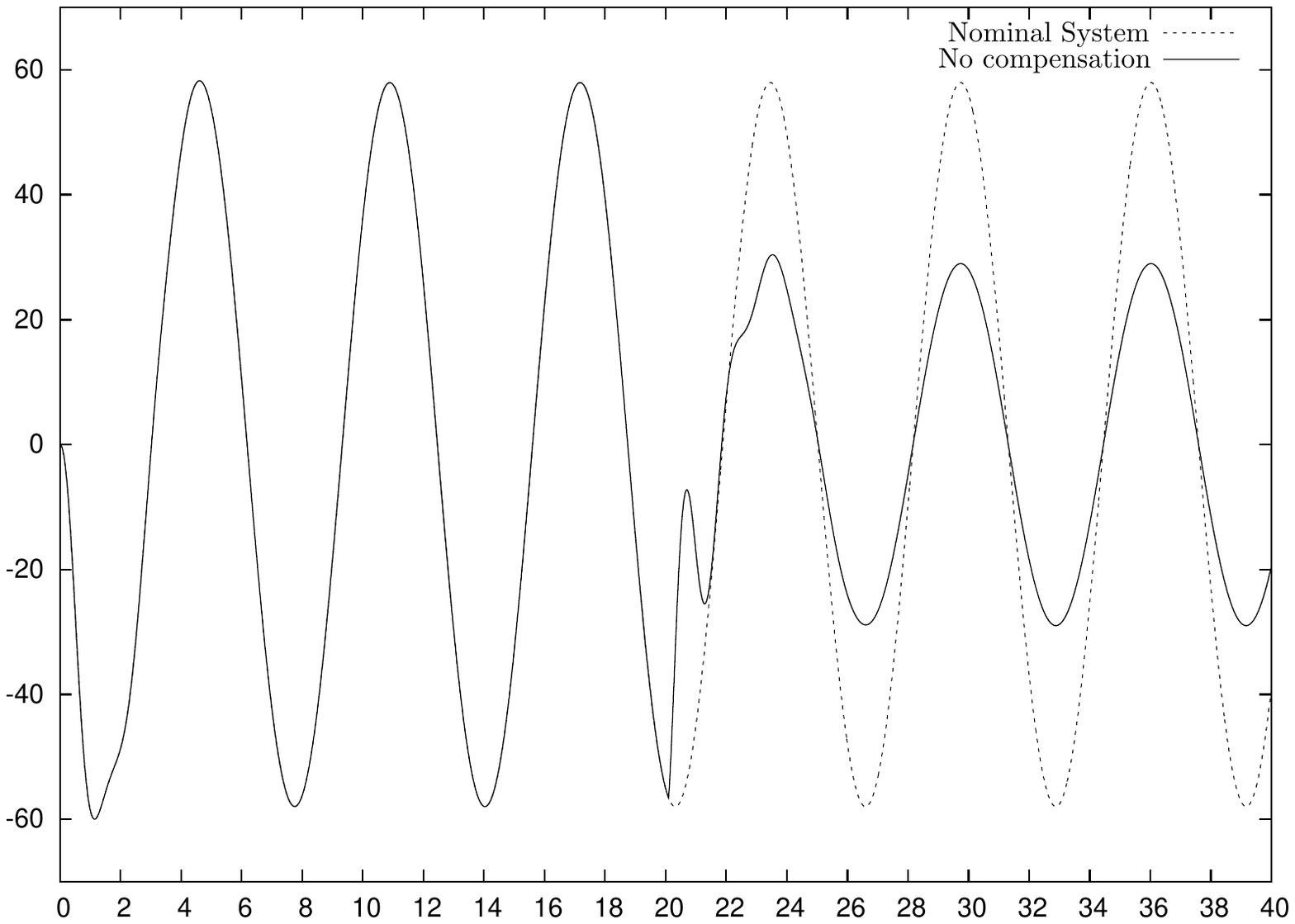}
\includegraphics[width=6cm,viewport=225 110 680 450]{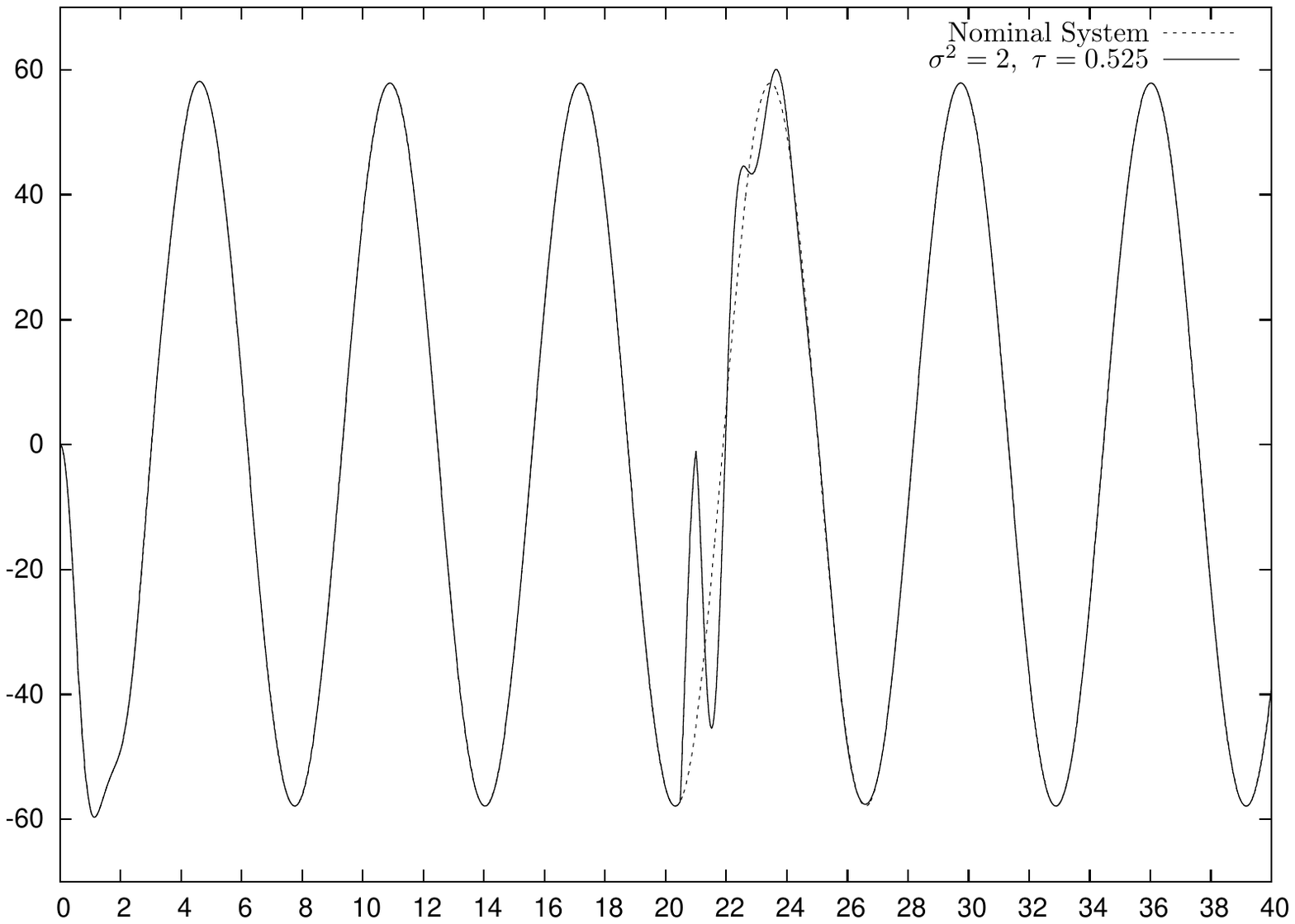}
\includegraphics[width=6cm, viewport=225 110 680 450]{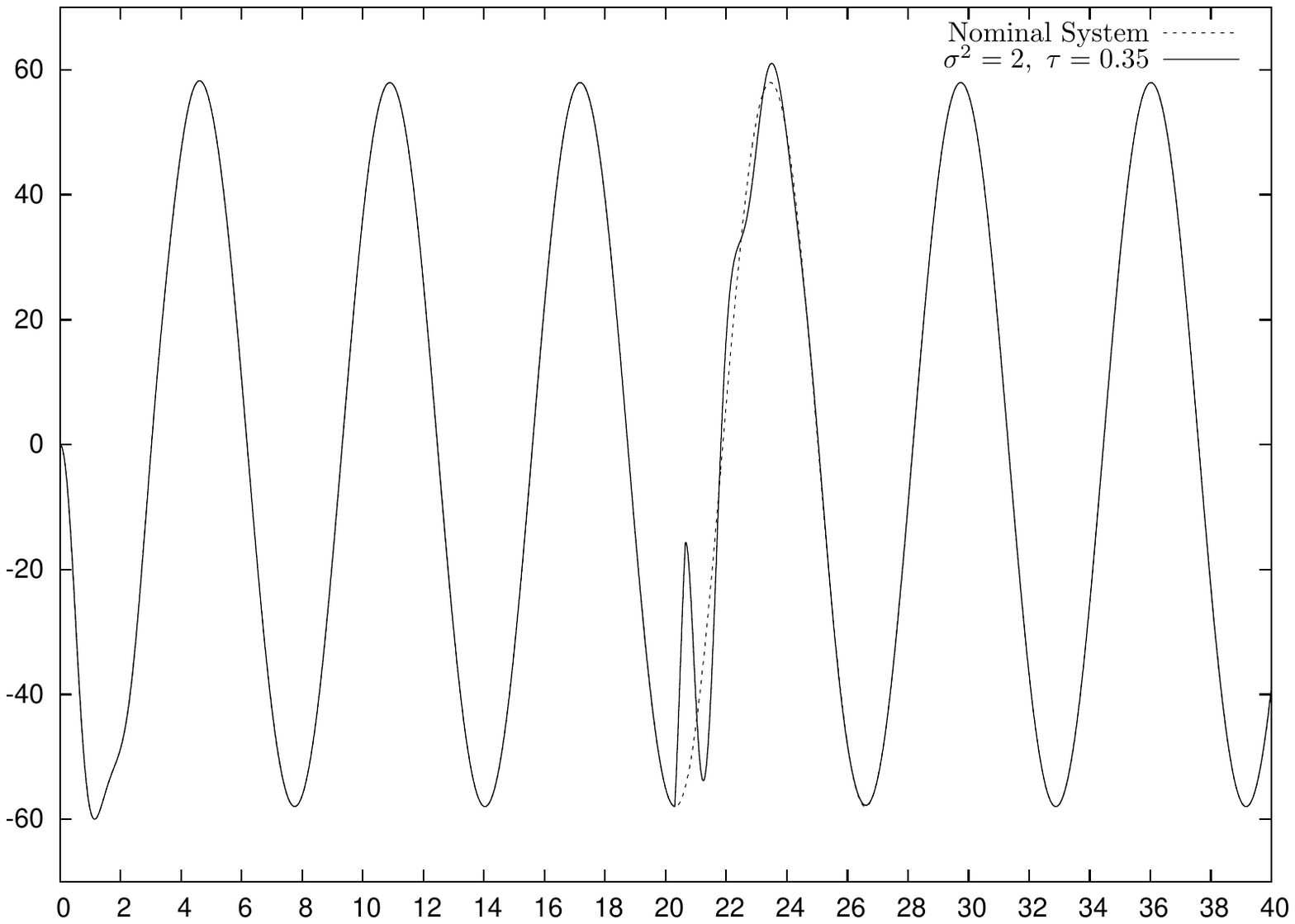}
\includegraphics[width=6cm,viewport=225 110 680 450]{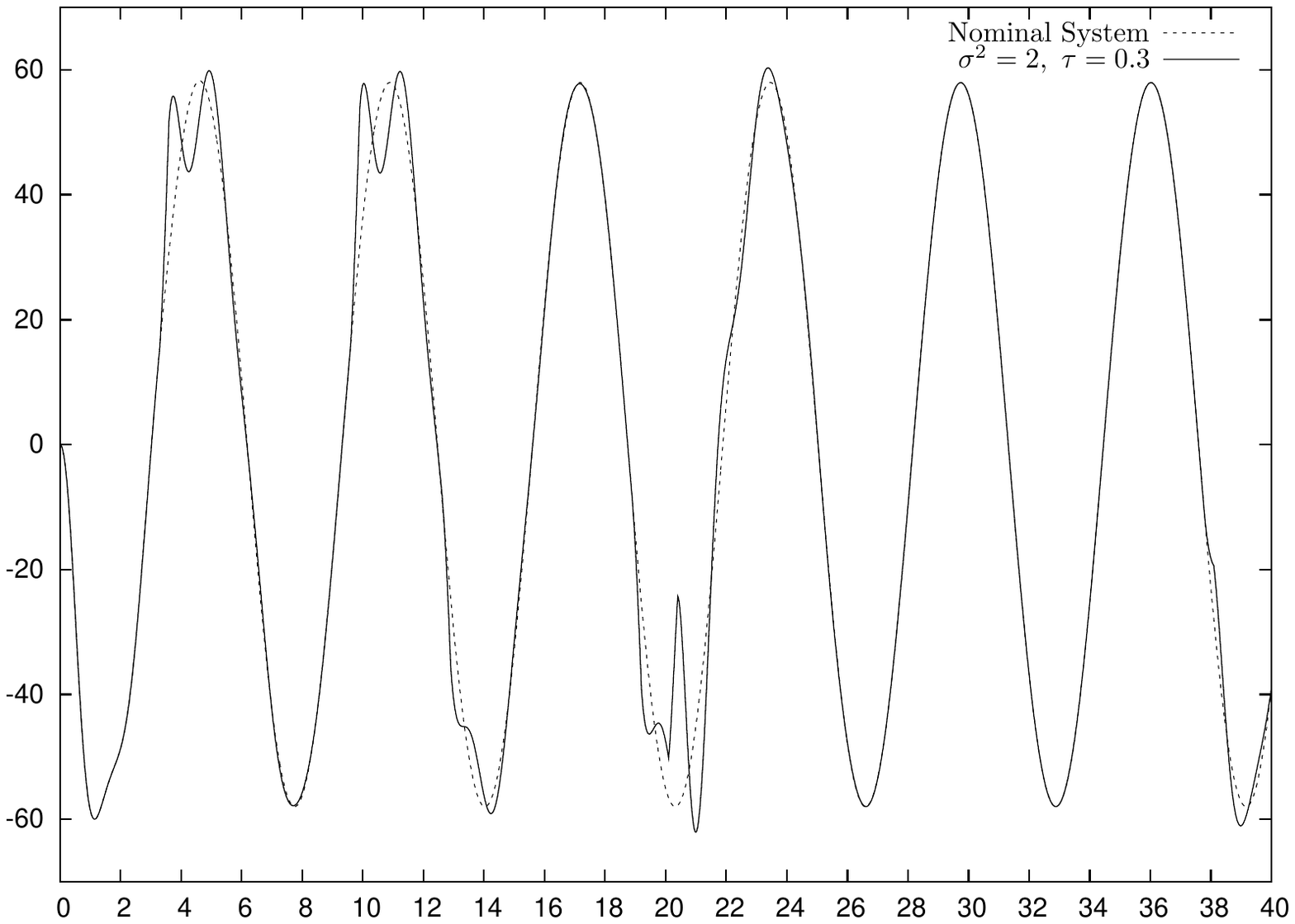}
\includegraphics[width=6cm,viewport=225 110 680 450]{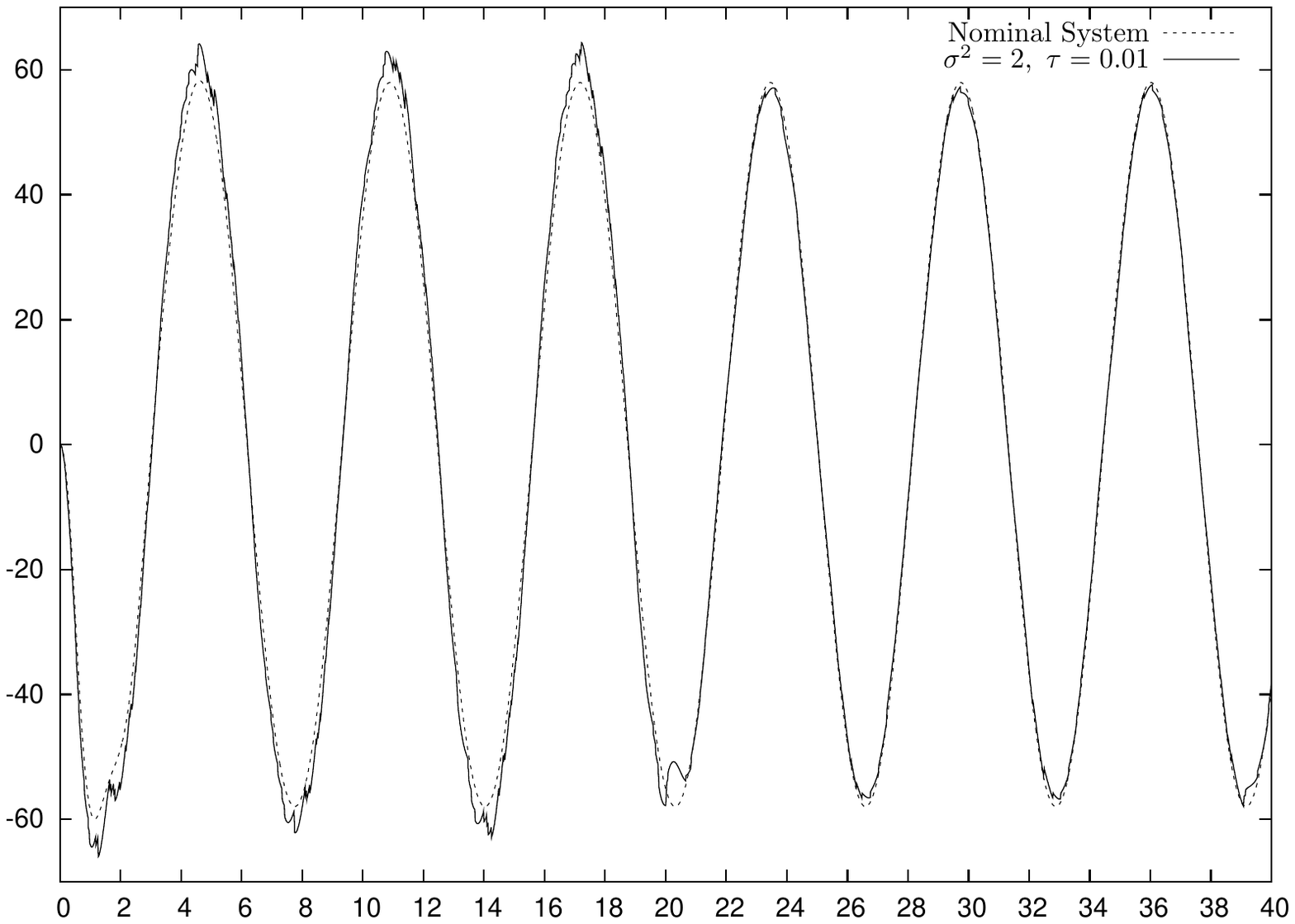}
\includegraphics[width=6cm,viewport=225 110 680 450]{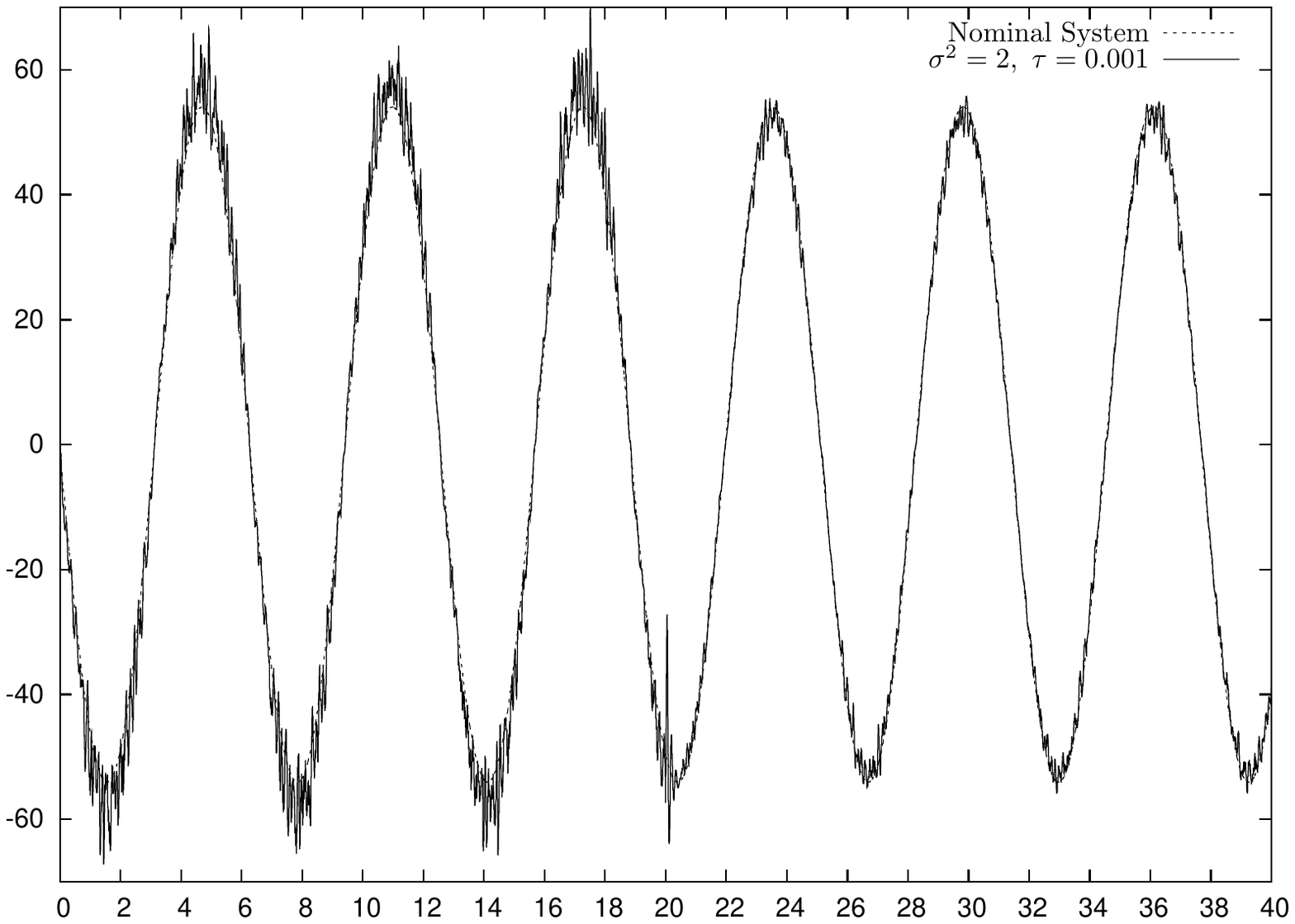}
\caption{Output y(t): Nominal System vs System with a failure at $T_F=20$, with lecture noise of variance $\sigma^2=2$ and $f(t)=\sin
t$. Six different cases are shown: the first graph represents the system with no control and compensation; the other ones are with
compensation, respectively with time step $\tau$ equal to $0.525$, $0.35$, $0.3$, $0.01$, $0.001$}\label{F2}
\end{figure}

Figure \ref{F2} concerns the case $f(t)=\sin t$. Again, the output of the system with no compensation in the first graph undergoes an
evident change after the failure at $T_F=20$. Instead, applying the One State Algorithm
with time step $\tau=0.525$ (this value being suggested by the numerical computation of the EDP) allows to recover the nominal condition. The same occurs with
$\tau=0.35$, which is preferable for the smaller amplitude of the unavoidable deviation in correspondence to the switch point.

When $\tau=0.3$, some detections fail (the error percentage is about $4\%$), but the output $y$ is not dramatically affected by them.
Furthermore, when $\tau=0.01$ the error percentage is about $9\%$: many deviations occur, but they are not very large. In particular,
they are quite null when the slope of $y(t)$ is steeper. In
correspondence to the switch point a plain  oscillation is present, but it is less remarkable than in the cases of larger $\tau$.

Decreasing $\tau$ again, the percentage of wrong detections does not overpass $10\%$, but for very small values of $\tau$,  the system
is  unstable (see for
instance, the last graph corresponding to $\tau=0.001$) and many oscillations occur.

\section{Conclusions}
In this paper, an original Fault Tolerant Control method, based on Information and Coding Theory, has been introduced. Given a
linear system with a disturbance and supposing that the disturbance function is quantized over two levels, the detection task can be
tackled by decoding techniques. In particular, we have introduced the One State Algorithm which is a low-complexity, recursive decoding
algorithm, derived from the BCJR. Its
application to a Flight FTC problem has generated satisfactory outcomes even in case of relative high noise in the data acquisition.\\ 

The low-complexity encourages the implementation of this method; moreover, adjusting the sampling time step $\tau$, one can improve
its performance, according to the different values of noise and of input $f$. In some cases, for instance when $f$ is constant, an
optimal value of $\tau $ can be analytically computed with sufficient precision, where the optimality is intended in terms of trade-off
between convergence conditions  and amplitude of the deviations. Other arrangements might be obtained changing the values and the number
of levels of quantization.

\bibliographystyle{plain}
\bibliography{my}

  \end{document}